%% file: mainfile_arxiv.tex
\documentclass[11pt]{article}
\usepackage[total={8.5in,11in}, top=1in, left=1in, right=1in, bottom=1in, text={6.5in,9in}]{geometry}
\pdfoutput=1

\usepackage{url}

\usepackage{algorithmic}
\usepackage[dvips]{graphicx}
\usepackage{epsfig}
\usepackage{amsmath, amsthm}
\usepackage{color}

\newtheorem{theorem}{\indent Theorem}[section]
\newtheorem{lemma}{\indent Lemma}[section]

\newtheorem{proposition}{\indent Proposition}[section]
\newtheorem{corollary}[theorem]{\indent Corollary}

\long\def\old#1{}

\def\red#1{{#1}}
\def\blue#1{{#1}}
\def\green#1{{#1}}
\def\blue#1{{#1}}

\def\aoe{}
\def\ao#1{{#1}}

\def\rt#1{#1}  
\def\mod#1{#1} 
\long\def\jnt#1{{#1}}

\def\jblue#1{{#1}}
\def\jred#1{{#1}}
\def\jmj#1{{#1}}

\definecolor{orange}{rgb}{1,0.5,0}
\long\def\ora#1{{#1}}
\def\aoc#1{{#1}}

\def\jh#1{{#1}}

\long\def\old#1{}

\def\red#1{{#1}}
\def\blue#1{{#1}}
\def\green#1{{#1}}
\def\blue#1{{#1}}

\def\aoe{}
\def\ao#1{{#1}}

\def\rt#1{#1}  
\def\mod#1{#1} 
\long\def\jnt#1{{#1}}

\def\jblue#1{{#1}}
\def\jred#1{{#1}}
\def\jmj#1{{#1}}

\definecolor{orange}{rgb}{1,0.5,0}
\long\def\ora#1{{#1}}
\def\aoc#1{{#1}}

\def\jntr#1{{#1}}
\providecommand{\rev[1]}{{#1}}
\providecommand{\aorev[1]}{{#1}}

\def\oM{\overline{M}}
\def\oT{\overline{T}}

\providecommand{\floor[1]}{\lfloor  #1 \rfloor }
\providecommand{\ceil[1]}{\lceil #1  \rceil}
\def\Rin{\mbox{Rin}}
\def\Rout{\mbox{Rout}}

\providecommand{\comjh[1]}{\com{JH}{#1}}
\providecommand{\com[2]}{\begin{tt}[#1: #2]\end{tt}}

\begin{document}

\title{Distributed anonymous discrete function computation \jnt{and averaging}\thanks{\rev{J.\ Hendrickx is with the Universit\'e catholique de Louvain, Louvain-la-Neuve, Belgium;  {\tt\small julien.hendrickx @uclouvain.be}.  \aorev{A. Olshevsky is with the Department of Mechanical
and Aerospace Engineering, Princeton University, Princeton, NJ, USA; {\tt \small  aolshevs@princeton.edu}.} J.\ N.\ Tsitsiklis is} with the Laboratory for Information
and Decision Systems, Massachusetts Institute of Technology, Cambridge, MA, USA{\tt\small jnt@mit.edu.}
This research
was supported by the National Science Foundation under
\jnt{a graduate fellowship and}
grant
ECCS-0701623, and \jnt{by}
postdoctoral fellowships from the F.R.S.-FNRS (Belgian Fund for
Scientific Research) and the B.A.E.F. (Belgian American Education
Foundation), \rev{and was conducted while J.\ Hendrickx \jntr{and A.\ Olshevsky were} at M.I.T.} 
\jnt{A preliminary version of this paper was presented at the Forty-Seventh Annual Allerton Conference on Communication, Control, and Computing, October 2009.}
}}

\author{Julien M. Hendrickx, Alex Olshevsky, \mod{and} John N. Tsitsiklis}
\maketitle

\input{abstract}
\input{intro}
\input{negative_results}
\input{positive_results}
\input{maxtracking}

\input{average_computation}
\input{simulations}
\input{conclusions}

\input{biblio}
\newpage
\appendix
\input{appen_proof_max}
\input{complexity}



\end{document}

%% file: abstract.tex
\begin{abstract}
\noindent
We propose a model for \red{deterministic} distributed function computation by a
network of identical \green{and} anonymous nodes. \jnt{In this model, each node has} bounded computation
and storage \green{capabilities} that do not \jblue{grow} with the network size.
\jnt{Furthermore, each node only knows its neighbors, not the entire
graph.}
Our goal is to \green{characterize the class of functions that can
be computed within this model.}  \aoc{In our main result, we provide a
necessary condition for computability which we show to be nearly sufficient, in the sense that
every function that violates this condition can at least be approximated.}  \blue{The problem
of computing \jnt{(suitably rounded)} averages in a distributed manner plays a central role
in our development;} \mod{we provide an algorithm that solves
it in time \jnt{that grows} quadratically with the size of the network.}
\end{abstract}

%% file: intro.tex
\section{Introduction}\label{s:intro}

The goal of many multi-agent systems,
\green{distributed} computation \green{algorithms,} and
decentralized data fusion methods is to have a
set of \green{nodes compute} \blue{a} common value based on
\green{initial values or observations at each node. Towards this purpose,
the nodes,} which we will \green{sometimes refer to as agents,} perform some internal
computations and repeatedly communicate with each other.
\jnt{The objective of this paper is to understand the fundamental limitations and capabilities of such systems and algorithms when the
available information and computational resources at each node are limited.}

\subsection{\jnt{Motivation}}\label{se:motiv}

\ora{The model that we will employ is a natural one for many different settings, including the case of wireless sensor networks.
However, before describing the model, we}
\jnt{start with a few examples that motivate the questions that we address.}

\smallskip
\noindent {\bf (a) Quantized consensus:} Suppose
\green{that} each \jnt{node} begins with an integer value $x_i(0)
\green{\in \{0,\ldots,K\}}$. We would like \green{the \jnt{nodes} to
end up, at some later time, with values $y_i$ that are almost
equal, i.e., $|y_i - y_j| \leq 1$, for all $i,j$, while preserving
the sum of the values, i.e., $\sum_{i=1}^n x_i(0) =
\sum_{\blue{i=1}}^n y_i$.} This is the so-called quantized
\red{averaging} problem, which has received considerable attention
recently; see, \jnt{e.g.,} \cite{KBS07,FCFZ08,ACR07, KM08, BTV08}.
It may be viewed as
the problem of computing the function $(1/n) \sum_{i=1}^n x_i$,
rounded to \jnt{an integer value}. 

\smallskip
\noindent {\bf (b) Distributed hypothesis testing} and \green{\bf majority voting:} Consider $n$
sensors \green{interested in deciding} between two hypotheses,
$H_0$ and $H_1$. Each sensor collects measurements and makes a
preliminary decision \red{$x_i\in\{0,1\}$} in favor of one of the
hypotheses. The sensors would like to \green{make a final decision
by majority vote, in which case they need to compute the indicator
function of the event $\sum_{i=1}^n x_i \geq n/2$, in a
distributed way. Alternatively, in a weighted majority vote, they
may be interested in computing the indicator function of \jnt{an event such as}
$\sum_{i=1}^n x_i \geq 3n/4$.} \rt{A variation of this problem
involves the possibility that some sensors \jnt{abstain from the vote, perhaps} due to their inability to gather
sufficiently reliable information.}

\smallskip
\noindent {\bf (c) \aoc{Direction coordination on a ring}:}
\rt{Consider $n$ vehicles placed on a ring, each with some
arbitrarily chosen 
direction of motion (clockwise or counterclockwise). We would
like the $n$ vehicles to agree on a single direction of motion.}
\rt{A variation of this problem was considered in \cite{MBCF07},
where, however, additional requirements on the vehicles were
imposed which we do not consider here. The solution \jnt{provided} in
\cite{MBCF07} was semi-centralized in the sense that vehicles had
unique \jnt{numerical} identifiers, and the final direction of most \jnt{vehicles} was \jnt{set to the
direction of the vehicle} with the \jnt{largest identifier.}  We wonder \jnt{whether} the
direction coordination problem \jnt{can} be solved in a completely
decentralized way. Furthermore, we would like the final direction
of motion to correspond to the initial direction of the majority
of the vehicles: if, say, 90\% of the vehicles are moving
counterclockwise, \jnt{we would like} the other 10\%  to turn around.}
\rt{If we define $x_i$ to be $1$ \jnt{when} the $i$th vehicle is
initially oriented clockwise, and $0$ if it is oriented
counterclockwise, then, coordinating on a direction involves the
distributed computation of the \jnt{indicator function of the event}  $\sum_{i=1}^n x_i \geq
n/2$}.

\smallskip
\noindent {\bf (d) Solitude verification:} \green{This is the
problem of} \jnt{checking whether exactly} \green{one} node has a given state. \green{This problem is of interest if we want
to avoid simultaneous transmissions over a common channel
\cite{GFL83}, or if we want to maintain a single leader (as in
motion coordination  --- see for example \cite{JLM03}) Given
\aoc{states} $x_i \in
\{\aoc{0},1,\ldots,\blue{K}\}$, \green{solitude verification} is
equivalent to \green{the problem of} computing the binary function
which is equal to 1 if and only if $|\{i: x_i=\aoc{0}\}| =1$.}
\smallskip

There are numerous methods that have been proposed for
solving problems such as the above; see for example the vast and
growing literature on consensus and averaging methods, \jntr{or the distribute robotics literature \cite{BCM09}.} 
Oftentimes, different algorithms involve different computational
capabilities on the part of the \jnt{nodes}, which makes it hard to
talk about \jnt{a} ``best'' algorithm. At the same time, simple
algorithms (such as setting up a spanning tree and \jnt{aggregating}
information by progressive summations over the tree) are often
\jnt{considered undesirable} \blue{because they require} too much coordination
or global information. \jnt{It should be clear} that a sound discussion
of such issues requires the specification of a precise model of
computation, followed by a systematic analysis of fundamental
limitations under  a given model. This is precisely the
objective of this paper: \blue{to} propose a particular model, and
to characterize the class of \jnt{functions computable} under this
model.

%
%


\subsection{\ora{The features of our model}}\label{se:features}

Our model provides an abstraction for \aoc{common}
requirements for distributed algorithms in the \ora{wireless} sensor network
literature.
\jnt{We
model the nodes as interacting
deterministic finite automata that
exchange messages on a fixed bidirectional network,
with no time delays or unreliable transmissions.
Some important qualitative features of our model are the following.}

\smallskip

\noindent \textbf{\red{Identical} \jnt{nodes}:} Any two \jnt{nodes} with the
same number of neighbors must run the same algorithm. \rev{Note that this assumption is equivalent to assuming that the nodes are exactly identical.
\aorev{Any algorithm \jntr{that} works in this setting will also work if the nodes are not all identical, since the nodes can still
run the same algorithm.}} 

\smallskip

\noindent\textbf{Anonymity:} \green{\jnt{A node} can distinguish its
neighbors using its own, private, local identifiers. However,
\jnt{nodes} do not have global identifiers.} \rev{In other words, a node receiving a message from one of its neighbors can send an answer to precisely that neighbor, or recogni\jntr{z}e that a later message comes from this same neighbor. On the other hand, nodes do not a priori have a unique signature that \jntr{can} be recogni\jntr{z}ed by every other node.}

\smallskip

\noindent\textbf{\jnt{Determinism:}}
\jnt{Randomization is not allowed. This restriction is imposed in order to preclude essentially centralized solutions that rely on randomly generated distinct identifiers and thus bypass the anonymity requirement. Clearly, developing algorithms is much harder, and sometimes impossible, when randomization is disallowed.}

\smallskip
\jnt{
\noindent\textbf{Limited memory:} We focus on the case where the nodes can be described by finite automata, and pay special attention to the required memory size. Ideally, the number of memory bits required at each node should be bounded above by a slowly growing function of the degree of a node.}

\smallskip

\noindent\textbf{Absence of global information:} \jnt{Nodes} have no
\green{global information, and do not even have an upper bound on
the total number of nodes. Accordingly, the algorithm that each
\jnt{node} is running is independent of the network size and
topology.}\smallskip

\noindent\textbf{Convergence \jnt{requirements:}} \jnt{Nodes hold an estimated output that
must converge to a desired value which is a function of all nodes' initial
observations or values.}
 \green{In particular, for the case
of discrete outputs,} all \jnt{nodes} must eventually \green{settle on}
the desired value. \green{On the other hand, the \jnt{nodes} do not
need to \jnt{become} aware of such termination, which is anyway impossible
in} the absence of \green{any} global information \cite{ASW88}.

\smallskip

\jnt{\jblue{In this paper, we
only consider the special case of} \textbf{fixed graph topologies,} where the underlying (and unknown) interconnection graph does not change with time. Developing a meaningful model for the time-varying case and extending our algorithms to that case is \jblue{an interesting topic, but} outside the scope of this paper.}

\subsection{\jnt{Literature review}}

There is a \jnt{very} large literature on distributed function
computation in related models of computation \cite{BT89, L96}.
\jnt{This literature can be broadly divided into two strands,
although the separation is not sharp: works that address general
computability issues for various  models, and works that focus on
the computation of specific functions, such as the majority function
or the average. We start by discussing the first strand.}

A common model in the distributed computing literature involves the
requirement that \green{all} processes terminate \jnt{once} the
\green{desired output is produced} \jblue{and \aoc{that} nodes become
aware
that termination has occurred.} A consequence of the termination
requirement is that nodes typically need to know the network size
$n$ (or an upper bound on $n$) to compute non-\aoc{trivial} functions. We
refer the reader to \cite{A80, ASW88, YK88, KKB90, MW93} for
\green{some} fundamental results in this setting, \green{and to
\cite{FR03} for} a \jnt{comprehensive} summary of known results.
\ora{Closest to our work is the reference \cite{CS}  which provides
an impossibility result very similar to our Theorem
\ref{unbounded-bound}, for a closely related model computation.}

\ora{The} biologically-inspired ``population algorithm" model
of distributed computation has some features in common with our
model, namely, \jblue{anonymous,} \jnt{bounded-resource} \jnt{nodes,} and
\jblue{no requirement of termination awareness;}
\green{see \cite{AR07} for \jnt{an overview of available}
results.} However, \green{\jnt{this} model involves a
different type of node interactions from the ones we consider;
in particular, \jblue{nodes interact pairwise at times that may be chosen adversarially.}
}

\jnt{Regarding the computation of specific functions, \cite{LB95}
shows} the impossibility of \jnt{majority voting} \jnt{if the nodes}
\jblue{are limited to a binary state.} Some experimental memoryless
algorithms \jnt{(which are not guaranteed to always converge to the
correct answer) have been} proposed \green{in the physics
literature} \cite{GKL78}.  Several papers have quantified the
performance of simple heuristics for computing specific functions,
typically in randomized settings. We refer the reader \blue{to}
\cite{HP01}, \aoc{which} studied simple heuristics \green{for computing}
the
majority function, \jblue{and to \cite{PVV08}, which provides a
heuristic that has guarantees only for the case of complete graphs.}

The large literature on quantized
averaging often tends to involve themes \jnt{similar to those addressed in this paper}
\cite{FCFZ08, KM08, ACR07, CB08, KBS07}.  \jnt{However, the underlying models of computation
 are typically more powerful than ours, as they allow for randomization and unbounded memory.
  Closer to the current paper, \cite{NOOT07}  develops \ora{an algorithm with $O(n^2)$ convergence time}
  for a variant of the quantized averaging problem, but requires unbounded memory.
  Reference \cite{BTV08} provides an algorithm for the particular quantized averaging problem that
   we consider in Section \ref{sec:reduction_to_avg} \jblue{(called in \cite{BTV08} the ``interval consensus problem''),} which uses
\ora{randomization but only} bounded memory (a total of two bits at
each node). 
\rev{An upper bound on its expected convergence time is provided in \cite{DV10} as a function of $n$ and a spectral quantity related to the network.
A precise convergence time bound, as a function of $n$, \aorev{is not given.}} Similarly, the algorithm in \cite{ZM09} runs in
$O(n^5)$ time for the case of fixed graphs. (However, we note that
\cite{ZM09} also addresses \jntr{an asynchronous model involving} time-varying graphs.) Roughly
speaking, the algorithms in \cite{BTV08, ZM09} work by having
positive and negative ``load tokens'' circulate randomly in the
network until they meet and annihilate each other. Our algorithm
involves a similar idea. However, at the cost of some algorithmic
complexity, our algorithm is deterministic. This allows for fast
progress, in contrast to the slow progress of algorithms that need
to wait until the coalescence time of two independent random walks.}
Finally, a deterministic algorithm for computing the majority
function (and some more general functions) was proposed in
\cite{LBRA04}. However, the algorithm appears to rely on
\aoc{the computation of shortest path lengths}, and thus \aoc{requires unbounded
memory at each node.}

Semi-centralized versions of the problem, in which the nodes
ultimately transmit to a fusion center, have often \jblue{been} considered
in the literature, e.g., for distributed statistical inference
\cite{MK08} \green{or detection} \cite{KLM08}. The papers
\cite{GK05},
 \cite{KKK05}, and \cite{YSG07} consider the complexity of computing a function and
communicating \green{its value} to a sink node. We refer the reader
to the references therein for an overview of \green{existing}
results in such semi-centralized settings. \jblue{However, the
underlying model is fundamentally different from ours, because the
presence of a fusion center violates our anonymity assumption.}

\jblue{Broadly speaking,} our results differ from previous works in
several key respects: (i) Our model, which involves totally
decentralized computation, deterministic algorithms, and constraints
on memory and computation resources at the nodes, but does not
require the nodes to know when the computation is over, \ora{is
different from that considered in \aoc{almost all} of the relevant
literature.}
 (ii) Our focus is on identifying computable
\green{and} non-computable functions \blue{under our model,} and
we achieve a nearly tight separation, \rt{as evidenced by \jblue{a comparison} between Theorem \ref{unbounded-bound} and Corollary
\ref{cor:approx_set_epsilon}.}
\jnt{(iii) Our $O(n^2)$ \aoc{averaging} algorithm is quite different,
and \aoc{significantly} faster than available memory-limited algorithms.}

\subsection{\jnt{Summary and Contributions}}

We \green{provide}  a general model of decentralized
anonymous computation \jblue{on fixed graphs,} \jnt{with the features described in Section \ref{se:features},}
and characterize the type of functions of the initial
\green{values} that can be computed.

We prove that if a function is computable under \green{our model,}
then its value \jnt{can only depend} on the \green{frequencies of the
different possible initial values.} For example, if the initial
\green{values} $x_i$ \jnt{are binary,} a computable
function \jnt{can} only depend on $p_0:= | \{i:x_i=0\}|/n$ and
$p_1:=|\{i:x_i=1\}|/n$. In particular, determining the number of
nodes, or \green{whether} at least two nodes have an initial
\green{value} \blue{of} $1$, is \green{impossible.}

Conversely, we prove that if a function only depends on the
\green{frequencies of the different possible initial values (and
is measurable)}, then the function can be approximated
\green{with any given} precision, except possibly on a set of
frequency vectors of  arbitrarily small volume. Moreover, if the
dependence on \green{these frequencies} can be expressed \jblue{through} a
combination of linear inequalities with rational coefficients,
then the function is computable exactly. In particular, the
functions involved in the quantized consensus, distributed
hypothesis testing, \rt{and direction coordination} examples are
computable, whereas the function involved in solitude verification
is not. Similarly, statistical measures such as the standard
deviation \jnt{of the distribution of the initial values} can be approximated with arbitrary
precision.
Finally, we show that with infinite memory, \red{the frequencies
of the different \jnt{initial} values} (i.e., $p_0$, $ p_1$ in the binary case) are
computable \jnt{exactly}, \aoc{thus obtaining a precise characterization of the
computable functions in this case.}

\jnt{The key to our positive results is a new algorithm for calculating the (\jblue{suitably} quantized) average of the initial values, which is of independent interest.
The algorithm does not involve randomization, requires only $O(n^2)$ time to terminate, and the memory (number of bits) required at each node is only logarithmic in the node's degree. In contrast, existing algorithms either require unbounded memory, or are significantly slower to converge.}

\subsection{\rt{Outline}}

In
Section \ref{sec:formal_descr}, we describe formally our model of computation. In Section
\ref{sec:nec_conditions}, we establish necessary conditions for a function
to be computable. In Section \ref{sec:reduction_to_avg}, \jnt{we provide}
sufficient conditions for a function to be
computable or approximable.
\jnt{Our positive results}
rely on an algorithm that keeps track of nodes
with maximal values, and an algorithm that calculates
\jnt{a suitably rounded average of the nodes' initial values; these} are described in
Sections \ref{sec:max_tracking} and \ref{sec:comput_average},
respectively.
\aoc{We provide some corroborating simulations in
Setion \ref{se:simul}, and
\jnt{we end with some} concluding remarks, in
Section \ref{sec:conclusions}}.

\section{Formal description of the model}\label{sec:formal_descr}

\jnt{Under our model, a distributed computing system consists of three elements:

\smallskip
\noindent {\bf (a) A network:} A network is a triple $(n,G,{\cal L})$,
where $n$ is the number of nodes, and $G=(V,E)$ is a \aoc{{\it connected bidirectional}} graph $G=(V,E)$ with $n$
nodes. (By bidirectional, we mean that \aoc{the graph is directed but}
if $(i,j) \in E$, then $(j,i) \in
E$.) We define $d(i)$ as the in-degree (and also out-degree, hence ``degree'' for short) of node $i$. Finally, $\cal L$ is a \emph{port labeling} which assigns a
\emph{port number} \ora{(a distinct integer in the set $\{\aoc{0,1,}\ldots,d(i)\}$)}  to each outgoing edge of
\jblue{any} node $i$.} \rev{Note that the unique \jntr{identifiers $i$ used to refer to nodes are only introduced for the purpose} of analysis, and are not part of the actual system. In particular, nodes do not know and cannot use their \jntr{identifiers.} }


\jnt{
\smallskip\noindent {\bf (b) Input and output sets:}
The input set is a finite set $X=\{0,1,\ldots,K\}$ to which the initial value of each node belongs. The output set is a finite set $Y$ to which the output of each node belongs.}

\jnt{
\smallskip\noindent {\bf (c) An algorithm:} An algorithm is defined as a family of finite
automata $(A_d)_{d = 1,2,\ldots}$, where the automaton $A_d$
describes the behavior of a node with degree $d$.}
The state of the automaton $A_d$ is \green{a} tuple $\jnt{[}x,z,y; (m_1,
\ldots, m_d)]$;  we will call $x \in X$ the
{\it initial value}, $z\in \red{Z_d}$ the {\it internal memory state},
$y\in Y$ the {\it output} or {\it estimated answer}, and $m_1, \ldots, m_d \in
M$ the \jnt{outgoing} {\it messages.} The sets $Z_d$ and $M$ are assumed finite.
\jnt{We allow the cardinality of $Z_d$ to increase with $d$. Clearly, this would be necessary for any algorithm
that needs to store the messages received in \aoc{the} previous time step.}
\jnt{Each automaton $A_d$ is identified with a transition law
from $X \times \blue{Z_d}
\times Y \times M^d$ into itself, which maps each}
$\left[x,z,y;(m_1,\dots,m_d)\right]$ to some
$\left[x, z', y';(m'_1,\dots, m'_d)\right]. $ In words, \jblue{at each iteration, the
automaton takes $x$, $z$, $y$, and incoming messages into account, to create} a new memory state, output, and
\green{(outgoing)} messages, but does not change
the initial \blue{value.}

\smallskip
\jnt{Given the above elements of a distributed computing system, an algorithm proceeds as follows.
For convenience, we assume that the above defined sets $Y$, $Z_d$, and
$M$ contain a special element, denoted by $\emptyset$.
Each} node $i$ begins with an initial value
$x_i\in X$ and implements the automaton $A_{d(i)}$, initialized with $x=x_i$ and
$z=y=m_1=\cdots=m_d = \emptyset$.
\green{We use $S_i(t)=
[ \blue{x_i},y_i(t),z_i(t),m_{i,1}(t), \ldots, m_{i, d(i)}(t)]$
to denote the}
state
of \jblue{node $i$'s} automaton at \jblue{time} $t$.
\jnt{Consider a particular node $i$. Let}
$j_1, \ldots,
j_{d(i)}$ be an enumeration of \jnt{its} neighbors, according to the port numbers. (Thus, $j_k$ is the \jblue{node at the other end of} the $k$th outgoing edge at node $i$.) Let $p_k$
be the port number \jnt{assigned to link $(j_k,i)$ according to the port labeling at node $j_k$. At each time step, node $i$ carries out the following update:}

\begin{small}
$$\left[x_i,
z_i(t+1),y_i(t+1); m_{i,1}(t+1), \ldots, m_{i, d(i)}(\red{t+1}) \right]  =
A_{d(i) }\left[\jnt{x_i}, z_i(t), y_i(t) ; m_{\blue{j_1, p_1}}(t),
\ldots, m_{j_{d\jnt{(i)}}, p_{d{\jnt{(i)}}}}(t) \right].
$$
\end{small}

\vspace{-8pt}
\noindent
In words, the
messages \jnt{$m_{j_{k}, p_{k}}(t)$, $k=1,\ldots,d(i)$,} ``sent'' by the  neighbors of $i$ into the ports
leading to $i$ are used to transition to a new state and
create new messages \jnt{$m_{i,k}(t+1)$, $k=1,\ldots,d(i)$,} that  $i$ ``sends'' to its neighbors at \jblue{time $t+1$.} We say that
\jnt{the algorithm \emph{terminates} if there exists some
$y^*\in Y$ (called the \emph{\green{final output}} of
the algorithm) and} a time $t'$ such that $y_i(t) =
y^*$ for every $i$ and $t\geq t'$.

Consider now a family of functions $(f_n)_{n = 1,2,\ldots}$, \jnt{where
$f_n:X^n\to Y$.} We say that such a
family is \emph{computable} if there exists a family of automata
$(A_d)_{d=1,2,\ldots}$ such that \blue{for any $n$,} for any
\jnt{network $(n,G,{\cal L})$},
and any set of initial conditions
$x_1,\dots,x_{\blue{n}}$, the \jnt{resulting algorithm terminates and the final output} is
$f_{\blue{n}}(x_1,\dots,x_n)$. \rev{Intuitively, a family of function\jntr{s} is computable \aorev{if there is a \jntr{bounded}-memory algorithm
which allows the nodes to ``eventually'' learn the value of the function in any connected topology and for any initial conditions.}}

\jnt{As an exception to the above definitions, we note that although we primarily focus on the finite case, we will} \jblue{briefly consider
in Section 
\rev{\ref{sec:reduction_to_avg}}}
function families $(f_n)_{n = 1,2,\ldots}$ {\em computable with infinite memory},  by which
we mean that the internal memory \blue{sets $Z_d$} and \jnt{the} output set $Y$ are
\jnt{countably infinite,} the rest of the model \jnt{remaining the same.}

\jnt{The rest of the paper focuses on the following general question:}
what \green{families} of functions are computable, and how can we
design a \jnt{corresponding algorithm $(A_d)_{d=1,2,\ldots}$?
To illustrate the nature of our model and the type of algorithms that it supports, we provide a simple example.}

\smallskip

\rt{\noindent \textbf{Detection problem:} In this \jnt{problem}, all
nodes \jnt{start with a binary initial value $x_i\in \{0,1\}=X$. We wish to
detect whether at least one node has an initial value equal to $1$.
We are thus dealing with the function family $(f_n)_{n=1,2,\ldots}$, where
$f_n(x_1,\ldots,x_n)=\max\{x_1,\ldots,x_n\}$.
This function family is computable by a family of
automata with binary messages, \jblue{binary} internal state, and
with the following transition rule:}}

\algsetup{indent=2em}
\begin{algorithmic}
\medskip
\IF {$x_i=1$ or $z_i(t) = 1$ or $\jnt{\max}_{j:(i,j)\in E}m_{ji}(t) =
1$}
 \STATE
set $z_i(t+1) = y_i(t+1) =1$ \STATE send $m_{ij}(t+1) = 1$ to
every neighbor $j$ of $i$

\ELSE \STATE set $z_i(t+1) = y_i(t+1) =0$ \STATE send $m_{ij}(t+1)
= 0$ to every neighbor $j$ of $i$

\ENDIF
\end{algorithmic}

\smallskip

\jnt{In the above algorithm, we initialize by setting $m_{ij}(0)$, $y_i(0)$, and $z_i(0)$ to zero instead of the special symbol $\emptyset$.}
One can easily verify that if $x_i=0$ for every $i$,
\jnt{then $y_i(t)=0$ for all $i$ and $t$.}
If on the other hand $x_k=1$
for some $k$, then at each time step $t$, those nodes \jnt{$i$} at distance
less than $t$ from \jnt{$k$} will have \jnt{$y_i(t)=1$.}
\jnt{Thus, for connected graphs, the algorithm will terminate within $n$ steps, with the correct output.}
It is important to note, \jnt{however, that because $n$ is unknown, a node $i$
can never know whether
its current output $y_i(t)$} is the final one. In particular, if
$y_i(t)=0$, node $i$
\jnt{cannot exclude the possibility that $x_k=1$ for some node whose distance from $i$ is larger than $t$.}

%% file: negative_results.tex
\section{Necessary condition for computability}\label{sec:nec_conditions}

\jnt{In this section we establish our main negative result, namely, that if a function family is computable, then the final output can only depend
on the frequencies of the different possible initial values. Furthermore, this remains true even if we allow for infinite memory,}
\ora{or restrict to networks in which neighboring nodes share a common label for the edges that join them. This result is quite similar to Theorem 3 of \cite{CS}, and so is the proof. Nevertheless, we provide a proof in order to keep the paper self-contained.}

\def\pb{frequency-based}

\jnt{We first need some definitions. Recall that $X=\{0,1,\ldots,K\}$. We let $D$ be the \emph{unit simplex,} that is,} $D =\{(\jnt{p_0},\dots,p_K)\in [0,1]^{\jnt{K+1}} : \sum_{k=\jnt{0}}^Kp_k =1 \}$. We say that a function $h:D\to Y$
\emph{corresponds} to a \jnt{function family $(f_n)_{n=1,2,\ldots}$ if for every $n$ and} every $x\in \blue{X}^n$, we have
$$
f(x_1,\dots,x_n)=
h\left(\aoc{p_0}(x_1,\ldots,x_n), p_1(x_1,\ldots,x_n),\dots,p_K(x_1,\ldots,x_n)\right),
$$
where
\[ p_k(x_1,\ldots,x_n) = |\{ \blue{i} ~|~ x_{\blue{i}} = \blue{k} \}|/n,\] so that
$p_k(x_1,\ldots,x_n)$ is the frequency of occurrence of the initial value $k$.
In this case, we say that the family $(f_n)$ is \emph{\pb.} \rev{Note that $n$ is used in defining the notion of frequency, but its value is unknown to the agents, and \jntr{cannot}  be used in the computations.}

\begin{theorem} \label{unbounded-bound}
Suppose that the family $(f_n)$ is computable with infinite memory.
Then, this family is \pb. \ora{The result remains true even if we
only require computability over edge-labeled networks.}
\end{theorem}

The following \jnt{are some applications of Theorem \ref{unbounded-bound}.}

\def\texitem#1{\par
\noindent\hangindent 25pt
\hbox to 25pt {\hss #1 ~}\ignorespaces}

\texitem{(a)}
The parity function $\sum_{i=1}^n x_i ~ ({\rm mod} \ k)$ is not
computable, for any $k>\aoc{1}$.
\texitem{(b)}
In a binary setting \blue{$(X=\{0,1\})$}, checking whether the number of nodes with
$x_i=1$ is \jnt{larger than or equal to} the number of nodes with $x_i=0$ \jnt{plus 10} is
not computable.
\texitem{(c)}
Solitude verification, i.e.,
checking whether $|i: \{x_i=0\}|=1$, is not computable.
\texitem{(d)}
An aggregate difference function such as
$\sum_{i<j} |x_i - x_j|$ is not computable, even \jnt{if it is to be} \green{calculated modulo $k$.}

\subsection{Proof of Theorem \ref{unbounded-bound}}

\ora{The proof of Theorem \ref{unbounded-bound} involves a particular degree-two network (a ring), in which all port numbers take values in the set $\{\aoc{0},1,2\}$, and in which any two
edges $(i,j)$ and $(j,i)$ have the same port number, \rev{as represented in Fig. \ref{rings}. More precisely, it \jntr{relies} on showing that two rings obtained by repeating, respectively, $k$ and $k'$ times the same sequences of nodes, edges, and initial conditions are algorithmically 
\jntr{i}ndistinguishable, and that any computable family of functions must thus take the same value on two such rings.}
The proof proceeds through}
\jnt{a sequence of intermediate results, starting
with} the following lemma, \rev{which essentially \jntr{reflects the fact that the node identifiers used} in our analysis do not influence the course of the algorithm.
It can be easily proved
by induction on time, and its proof is omitted. The second lemma states that the value taken by  computable (families of) functions may not depend on which particular node has which initial condition.}

\begin{lemma} \label{switching}
Suppose that $G=(\{1,\ldots,n\},E)$ and $G'=(\{1,\ldots,n\},E')$
are isomorphic; that is, there exists a permutation $\pi$ such
that $(i,j) \in E$ if and only if $(\pi(i),\pi(j)) \in E'$.
Furthermore, suppose that the \green{port label at node $i$ for the edge leading to} $j$ in $G$ is the
same as the \green{port label at node $\pi(i)$ for the edge leading to} $\pi(j)$ in $G'$.  Then, \green{the state
$S_i(t)$ resulting} from \blue{the initial values $x_1,\ldots,x_n$} \jnt{on}
the graph $G$ is the same as \blue{the state} $S_{\pi(i)}(t)$
\blue{resulting from the initial values $x_{\jnt{\pi^{-1}}(1)},\ldots,x_{\jnt{\pi^{-1}}(n)}$}
 \jnt{on} the graph $G'$.
\end{lemma}

\begin{lemma} \label{symmetry}
Suppose that the family \jnt{$(f_n)_{n=1,2,\ldots}$} is computable with infinite memory
\aoc{on edge-labeled networks}. Then, each $f_i$
is invariant under permutations of its arguments.
\end{lemma}
\begin{proof}
\blue{Let ${\pi_{ij}}$ be \jnt{the} permutation that swaps $i$ with $j$ \jnt{(leaving the other nodes intact)}; with a slight abuse of notation, we also denote by $\pi_{ij}$} the mapping from \blue{$X^n$ to $X^n$} that swaps the
$i$th and $j$th elements of a vector.
\jnt{(Note that $\pi_{ij}^{-1}=\pi_{ij}$.)}
We show that for all $x \in
\blue{X^n}$, $f_{\blue{n}}(x) = f_{\blue{n}}( \pi_{ij}( x))$.

We run our \green{distributed} algorithm \green{on} the $\blue{n}$-node
\aoc{complete graph with an edge labeling.} \aoc{Note that at least one
edge labeling for the complete graph exists: for example, nodes $i$ and $j$ can use port number $(i+j) ~{\rm mod }~ n$ for the edge connecting them.}. Consider two different \jnt{sets of} initial
\jnt{values, namely the vectors  (i) $x$,  and (ii) $\pi_{ij}(x)$.
Let the port labeling in} case (i) be arbitrary;
in case (ii), let the \jnt{port labeling be}
\blue{such that the conditions in Lemma \ref{switching} are satisfied  (which is easily accomplished).}
Since the \jnt{final value} is $f(x)$ \jnt{in case (i)} and
$f(\pi_{ij}(x))$ \jnt{in case (ii)}, we obtain $f(x)=f(\pi_{ij}(x))$.
\blue{Since the permutations $\pi_{ij}$ generate the group of permutations, permutation invariance follows.}
\end{proof}

Let \blue{$x\in X^n$.}  We will \green{denote by} $x^2$ the
concatenation of $x$ with itself, and, generally, \jnt{by} $x^k$ the
concatenation of $k$ copies of $x$. We now prove that
self-concatenation does not affect the value of a computable
family of functions.

\begin{lemma} \label{repetition} Suppose that the family
\jnt{$(f_n)_{n=1,2,\ldots}$} is computable with
infinite memory \aoc{on edge-labeled networks}. Then, \red{for every $n\geq 2$}, every sequence
\blue{$x\in X^n$}, and \green{every positive} integer~$m$, \[ f_n(x) = f_{mn} (x^m).\]
\end{lemma}

\begin{proof}
Consider a ring of \jnt{$n$ nodes,} where the $i$th \jnt{node}
\aoc{clockwise} begins with the $i$th element of $x$; and consider
a ring of \jnt{$mn$ nodes}, where the \jnt{nodes}
$i,i+n,i+2n,\ldots$ (\aoc{clockwise}) begin with the $i$th element
of $x$. Suppose that the labels in the first ring are $0,1,2,1,2,\ldots$.
That is, the label of the edge $(1,2)$ \jntr{is 0 at both nodes 1 and 2; the label of a subsequent edge $(i,i+1)$ is the same at 
both nodes $i$ and $i+1$, and alternates between 1 and 2 as $i$ increases.}  In the second ring, we simply repeat $m$ times the labels in
the first ring. See Figure \ref{rings} for an example with
$n=5$, $m=2$.

\begin{center}
\begin{figure}[h] \hspace{3cm}
\begin{minipage}[t]{.45\textwidth}
  \begin{center}        \epsfig{file=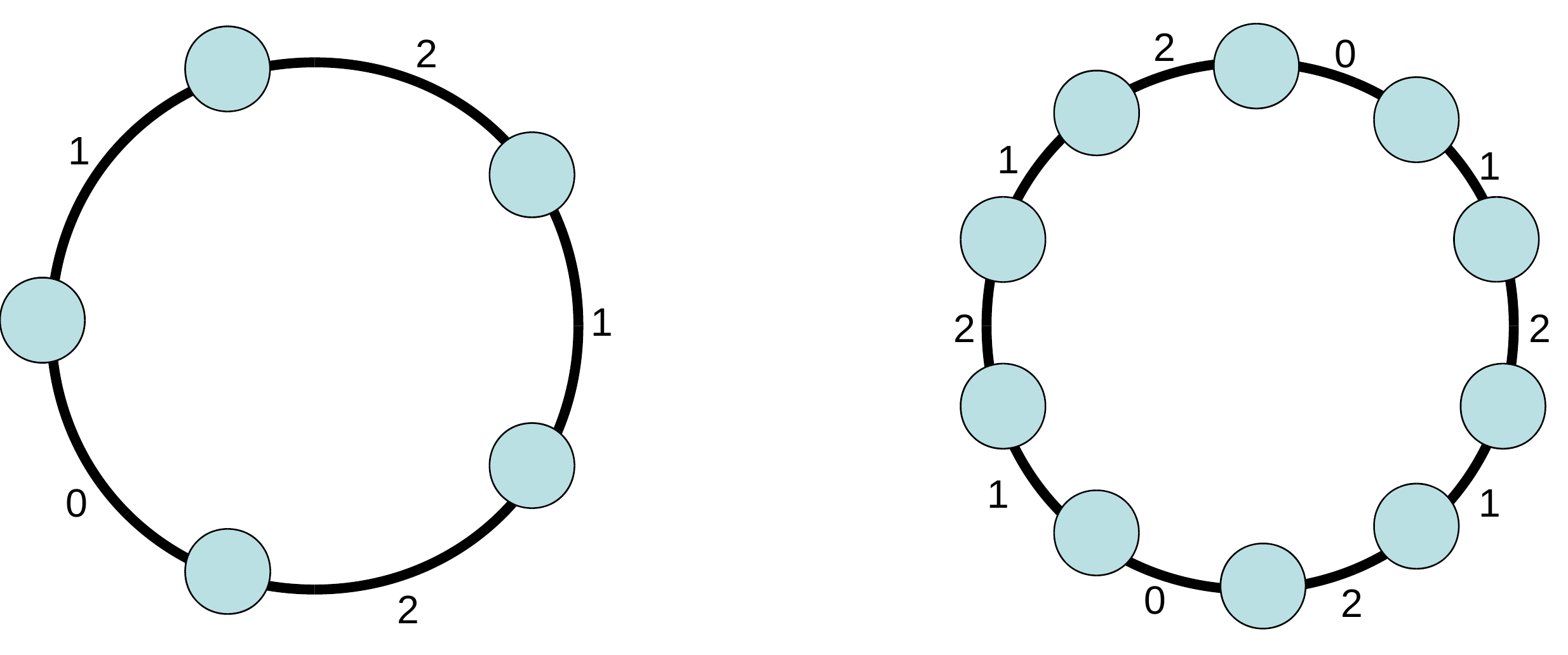, scale=0.45}
  \end{center}
\end{minipage}
\caption{Example of two situations that are \jnt{algorithmically} \blue{indistinguishable}. \ora{The numbers next to each edge are the edge labels.}
} \label{rings}
\end{figure}   \end{center}

Initially, the state \jnt{$S_i(t)=
[ \blue{x_i},y_i(t),z_i(t),m_{i,1}(t), m_{i, 2}(t)]$, with $t=0$,} of node $i$ in the first ring is exactly the
same as the state of the nodes \blue{$j=i,i+n,i+2n,\ldots$} in the second ring. We show by induction \green{that this property} must hold at all
times $t$. \blue{(To keep notation simple, we assume, without any real loss of generality,
 that $i\neq 1$ and $i\neq \aoc{n}$.)}

Indeed, suppose this \green{property} holds up to time $t$. At time $t$, node $i$ in
the first ring receives a message from node $i-1$ and a message
from node $i+1$; and in the second ring, node \red{$j$} satisfying $\red{j \ (\mbox{mod } n)  } =
i$ receives one message from $j-1$ and $j+1$.
Since $\red{j-1 \ (\mbox{mod } n)  } =
i-1 \aoe{\ (\mbox{mod } n) }$ and $\red{j+1 \ (\mbox{mod } n)  } =
i+1 \aoe{\ (\mbox{mod } n)}$,
the states of \green{$j-1$ and $i-1$} are identical at time $t$,
\green{and similarly for}
$j+1$ and $i+1$. Thus, \aoc{because of periodicity of the edge labels, nodes} $i$ (in the first ring) and $j$ (in the second ring) \aoe{receive
identical messages
\aoc{through identically labeled ports}} at time $t$. Since $i$ and $j$ were in the same state at time $t$, they
must be in the same state at time $t+1$. This proves that they are \jnt{always}
in the same state.
It follows that $y_i(t) = y_j(t)$ for all $t$, whenever
$\red{j \ (\mbox{mod } n)  } =
i$, and therefore $f_n(x)=f_{mn}(x^m)$.
\end{proof}

\begin{proof}[Proof of Theorem \ref{unbounded-bound}]
Let $x$ and $y$ be two sequences of $n$ and $m$
elements, respectively,  such that $p_k(x_1,\dots,x_n)$ \jnt{and} $p_k(y_1,\dots,y_m)$ \jnt{are equal to a common value} $\hat
p_k$, for $\jnt{k\in X}$; \jnt{thus,} the number of occurrences
of $k$ in $x$ and $y$ are $n\hat p_k$ and $m\hat
p_k$, respectively. Observe that \blue{for any $k\in X$, the vectors}  $x^m$ and $y^n$ have the same number $mn$ of
elements, and both contain $mn\hat p_k$ occurrences of \green{$k$.}
The sequences $y^n$ and $x^m$ can thus be obtained
from each other by a permutation, which by \jnt{Lemma}
\ref{symmetry} implies that $f_{nm}(x^m) = f_{nm}(y^n)$. From \jnt{Lemma} \ref{repetition}, \jnt{we have} that $f_{nm}(x^m) = f_n(x)$
and $f_{mn}(y^n)=f_m(y)$. \jnt{Therefore,} $f_n(x) = f_m(y)$. This proves
that the value of $f_n(x)$ is determined \aoc{by the values of $p_k(x_1,\ldots,x_n), k = 0,1, \ldots, K$.}
\end{proof}

\aorev{\noindent {\bf Remark:} Observe that the above proof remains valid even under the ``symmetr\jntr{y}'' assumption that 
\jntr{an edge is assigned the same label by both of the nodes that it is incident on.}
port label of each edge is the same at every node.}

%% file: positive_results.tex

\section{Reduction of generic functions to the computation of averages}\label{sec:reduction_to_avg}

\jnt{In this section, we turn to positive results, aiming at a converse of
Theorem
\ref{unbounded-bound}.
The centerpiece of our development is Theorem \ref{thm:average_computation}, which states that a certain average-like function is computable. Theorem \ref{thm:average_computation} then implies the computability of a large class of functions, yielding an approximate converse to Theorem
\ref{unbounded-bound}.} \rev{We will then illustrate these positive results on some examples.}

\jnt{The average-like functions that we consider correspond to the ``interval consensus'' problem studied in \cite{BTV08}. They are defined as follows.
Let $X = \{0,\ldots,K\}$. Let $Y$ be the following set of single-point sets and intervals:
\[ Y = \{ \{0\}, (0,1), \{1\}, (1,2),  \ldots,
\{K-1\}, (K-1,K), \{K\} \}\]
(or equivalently, an indexing of this finite collection of intervals).
For any $n$, let $f_n$ be the function that maps $(x_1,x_2,\ldots,x_n)$ to the element
of $Y$ which contains the average $\sum_i x_i/n$. We refer to the function family $(f_n)_{n=1,2,\ldots}$ as the \emph{interval-averaging} family.} \rev{The output of this family of functions is thus the exact average when it is an integer; \jntr{otherwise, it is} the open interval between two integers that contains the average. }

\jnt{The motivation for this function family comes from the fact that the exact average $\sum_i x_i/n$ takes values in a countably infinite set, and cannot be computed when the set $Y$ is finite. In the quantized averaging problem considered in the literature, one settles for an approximation of the average. However, such approximations do not necessarily define a single-valued function from $X^n$ into $Y$. In contrast, the above defined function $f_n$ is both single-valued and delivers an approximation with an error of size at most one. Note also that once the interval-average is computed, we can readily determine the value of the average rounded down to an integer.}

\begin{theorem}\label{thm:average_computation} \jnt{The interval-averaging function family}  is computable.
\end{theorem}


\jnt{The proof of Theorem \ref{thm:average_computation} (and the
corresponding algorithm) is quite involved; \ora{it will be
developed in Sections \ref{sec:max_tracking} and
\ref{sec:comput_average}.}} \aoc{In this section, we show that the
computation of a broad class of functions can be reduced to
interval-averaging.}

\jnt{Since only \pb\ function families can be computable (Theorem \ref{unbounded-bound}), we can restrict attention to the corresponding
functions $h$.}
\green{We
will say that a function $h$ on the \jnt{unit simplex $D$} is
\emph{computable} if it corresponds to a \pb\
computable family $(f_n)$.  The level sets of $h$ are defined as
the sets $L(y) = \{ \red{p} \in D ~|~ h(p) = y\}$, for $y \in Y$.}

\begin{theorem}[Sufficient condition for computability] \label{bounded-alg}
Let $h$ be a function from the \jnt{unit simplex} $D$ to $Y$. \green{Suppose that} every level set $L(y)$ can be written as a
finite union,
\[ L(y) =  \bigcup_k C_{i,k},\] where each $C_{i,k}$ can in turn be
written as a finite intersection of \green{linear} inequalities \green{of the form}
 \begin{equation} \aoc{\alpha_0 p_0} + \alpha_1 p_1 + \alpha_2 p_2 + \cdots + \alpha_K p_K \leq \alpha, \label{eq:ineq_first_type}
\end{equation} or
\[ \aoc{\alpha_0 p_0} + \alpha_1 p_1 + \alpha_2 p_2 + \cdots + \alpha_K p_K < \alpha,
\]

with rational \green{coefficients}  $\alpha,\aoc{\alpha_0}, \alpha_1,\ldots,\alpha_K$. Then, $h$
is computable.
\end{theorem}

\begin{proof}
Consider one such linear inequality, which we assume, for \aoc{concreteness}, to be of \jntr{the} \rev{form (\ref{eq:ineq_first_type})}. Let $P$ be the set of indices
\jnt{$k$} for which $\alpha_k \geq 0$. Since all coefficients are rational, we
can clear \aoc{their} denominators and rewrite \green{the inequality} as
\begin{equation}\label{eq:integ_equiv_density_ineq}
\sum_{k \in P} \beta_k p_k - \sum_{\jnt{k} \in P^c} \beta_k p_k \leq
\beta , \end{equation} for \red{nonnegative} integers $\beta_k$ and
$\beta$. Let $\chi_k$ be the indicator function \green{associated with initial
value} $k$, i.e., $\chi_k(i)=1$ if $x_i = k$, and $\chi_k(i) =0$
otherwise, so that $p_k=\frac{1}{n}\sum_{i}\chi_k(i)$.
\green{Then, (\ref{eq:integ_equiv_density_ineq}) becomes} $$\frac{1}{n} \sum_{i=1}^n \left(\sum_{k\in P}\beta_k
\chi_k(i) + \sum_{k\in P^c}\beta_k  (1-\chi_k(i))\right) \leq
\beta + \sum_{k\in \red{P^c}}\beta_k, $$ or \[ \frac{1}{n}\sum_{i\aoc{=1}}^{\aoc{n}} q_i
\leq  q^*, \] where $q_i = \sum_{k\in P}\beta_k \chi_k(i) +
\sum_{k\in P^c}\beta_k (1-\chi_k(i))$ and $q^* =\beta + \sum_{k\in
P^c}\beta_k$.

\rev{An algorithm \jntr{that determines}  whether the last inequality is satisfied can be designed \jntr{as follows.} Knowing the parameters $\beta_k$ and the set $P$, which can be made part of algorithm description as they depend on the problem and not on the data, each node can initially compute its value $q_i$, \jntr{as well as the value of $q^*$.} Nodes can then apply the distributed algorithm that
computes \jnt{the integer part of} $\frac{1}{n}\sum_{i=1}^n q_i$;
\jnt{this is possible by virtue of Theorem
\ref{thm:average_computation}, with $K$ set to $\sum_k\beta_k$ (the largest possible value of $q_i$).} It suffices then for them to  constantly compare \jntr{the current output} of this algorithm to the integer $q^*$.}   	
To check any finite collection of
inequalities, the nodes can perform the computations for each
inequality in parallel.

To compute $h$, the nodes simply need to check which \green{set} $L(y)$ the
\green{frequencies} $\aoc{p_0}, p_1,\ldots,p_K$ lie in, and this can be done by
checking the inequalities defining each $L(y)$.
\green{All of these computations} can be accomplished
with finite automata: indeed, we do nothing more than run finitely
many copies of the automata provided by Theorem
\ref{thm:average_computation}, one for each inequality.
\jnt{The total memory used by the automata depends on the number of sets $C_{i,k}$ and the magnitude of the coefficients $\beta_k$, but not on
$n$, as required.}
\end{proof}

Theorem \ref{bounded-alg} shows the computability of functions $h$
whose level-sets \green{can} be defined by linear inequalities
with rational coefficients. \red{On the other hand, it is clear
that not every function $h$ can be computable.} \blue{(This can be
shown by a counting argument: there are uncountably many possible
functions $h$ \jnt{on the rational elements of $D$,} but for the special case of bounded degree graphs,
only countably \aoc{many} possible algorithms.)} \green{Still, the next \jnt{result}
shows that the set of computable functions is rich enough, in the
sense that \jnt{computable}} functions can approximate any measurable
function, \jnt{everywhere except possibly on a low-volume set.}

We will call a set of the form $\prod_{\red{k}=\jnt{0}}^K (a_k,b_k)$,
\green{with every $a_k,b_k$ rational,} {\em a \red{rational open}
box}, where $\prod$ \green{stands for Cartesian} product. A
function that can be written as a finite sum $\sum_i a_i
\red{1_{B_i}}$, where \green{the $B_i$ are rational open boxes and
the $1_{B_i}$ are the associated indicator functions,} will be
referred to as a \emph{box function.} \green{Note that box
functions are computable by Theorem \ref{bounded-alg}.}

\begin{corollary}\label{cor:approx_set_epsilon}
If every level set of a function $h:D\to Y$ on the \jnt{unit simplex $D$}
is Lebesgue measurable, then, for every $\epsilon >0$, there exists a
computable box function $h_\epsilon:D\to Y$ such that the
set \aoc{$\{ p \in D ~|~ h(p) \neq h_{\epsilon}(p) \}$} has measure at most $\epsilon$.
\end{corollary}
\begin{proof} \jnt{(Outline)}
The proof relies on the following elementary result from measure
theory.
Given a Lebesgue measurable set $E \jnt{\subseteq D}$ and some $\epsilon
> 0$, there \green{exists a set $E'$ which is a finite union of
disjoint open boxes, and which satisfies}
\[ \mu( \aoc{E \Delta E'} ) < \epsilon,\] where $\mu$ is
the Lebesgue measure \rev{and $\Delta$ is the symmetric difference operator.}
\red{By a routine argument, these boxes can be taken to be rational.
By applying this fact to the level sets of the function $h$ (assumed measurable), the function $h$ can be approximated by a box function}  \blue{$h_{\epsilon}$. Since box functions are computable, the result follows.}
\end{proof}

The following corollary \green{states} that  
continuous functions are approximable.

\begin{corollary}\label{cor:approx_continuous}
If a function $h:D\to [L,U]\subseteq \Re$ is continuous, then for
every $\epsilon>0$ there exists a computable function
$h_\epsilon:D \blue{\to  [L,U]}$ such that
$\|h-h_\epsilon\|_\infty < \epsilon$
\end{corollary}
\begin{proof}
Since $D$ is compact, \jnt{$h$} is uniformly continuous. One can
therefore \green{partition} $D$ into a finite number of subsets,
$A_1,A_2,\dots,A_q$, that can be \green{described} by linear inequalities
with rational coefficients, so that $\max_{p\in
A_j}h(p)-\min_{p\in A_j}h(p) <\epsilon$ holds for all $A_j$. The
function $h_\epsilon$ is then built by assigning to each $A_j$ an
\green{appropriate} value \green{in} $\{L,L+\epsilon, L+2\epsilon, \dots, U\}$.
\end{proof}

To illustrate \jnt{these} results, let us consider again
some examples.

\texitem{(a)} \green{Majority \jnt{voting} between two options is} equivalent
to checking whether $p_1 \leq 1/2$, with alphabet $\{0,1\}$. \rev{This condition is clearly of the form (\ref{eq:ineq_first_type})},  and is
therefore computable.
\texitem{(b)}
Majority \jnt{voting} when some nodes can ``abstain'' amounts to checking whether $p_1 -
p_{\aoc{0}} \geq 0$, with \jnt{input set} $\jnt{X=}\{0,1,{\rm abstain}\}$. This function family is
computable.
\texitem{(c)}
We can ask for the second
most popular value out of four, for example. In this case, the
sets $A_i$ can be decomposed into constituent sets defined by
inequalities such as $p_2 \leq p_3 \leq p_4 \leq p_1$, each of which
obviously has rational coefficients. \rev{The level sets of the function can thus clearly be expressed as unions of sets defined by \jntr{a collection} of linear inequalities of the type (\ref{eq:ineq_first_type})}, so that the function is computable.
\texitem{(d)}
\green{For any subsets $I,I'$ of
$\{\aoc{0}, 1,\ldots,K\}$, the indicator function of the set
where} $\sum_{i \in I} p_i
\jnt{>} \sum_{i \in I'} p_i$  is computable. This is equivalent to checking
whether more nodes have a value in $I$ than do in $I'$.
\texitem{(e)}
The \green{indicator functions of the sets defined by} $p_1^2 \leq 1/2$ and
$p_1 \leq \pi/4$ are measurable, so they are approximable. We are
unable to say whether they are computable.
\texitem{(f)}
The \green{indicator function of the set defined by} $p_1 p_2 \leq 1/8$ is approximable, but we are unable
to say whether it is computable.

Finally, we show that with infinite memory, it is possible to
recover the exact \green{frequencies $p_k$.} \red{(Note that this
is impossible with finite memory, because $n$ is unbounded, and
the number of bits needed to represent $p_k$ is also unbounded.)}
\jnt{The main difficulty is that $p_k$ is a rational number whose denominator can be arbitrarily large, depending on the unknown value of $n$. The idea is to run separate algorithms for each possible value of the denominator (which is possible with infinite memory), and reconcile their results.}

\begin{theorem} \label{thm:infinite_mem_positive}
The vector $(\aoc{p_0}, p_1,\ldots,p_K)$ is computable with infinite memory.
\end{theorem}
\begin{proof}
We show that $p_1$ is computable exactly, which is sufficient to
prove the theorem. Consider the following algorithm, \jnt{to be referred to as $Q_m$,} parametrized
by a \green{positive integer} \red{$m$.} The \jnt{input} set $X_m$ \jnt{is} $\{0,1,\ldots,
m\}$ and the output set $Y_m$
\jnt{is the same as in the interval-averaging problem:}
$Y_m = \{ \{0\}, (0,1), \{1\}, (1,2), \{2\}, (2, 3),
\ldots, \{m-1\}, (m-1,m), \{m\} \}$. If $x_i=1$, then node sets
its initial value $x_{i,m}$ to $m$; else, the node sets its
initial value $x_{i,m}$ to $0$. \green{The algorithm computes} the
function family $(f_n)$ which maps $X_m^n$ to the element of $Y_m$
containing $(1/n) \sum_{i=1}^n x_{i,m}$, \green{which is possible, by} Theorem
\ref{thm:average_computation}.

The nodes run \green{the algorithms} $Q_m$ for every
\green{positive integer value of $m$, in an interleaved manner.
Namely, at each time step, a node runs one step of a particular
algorithm $Q_m$, according to} the following order:
\[ Q_1, \hspace{0.2cm} Q_1, Q_2, \hspace{0.2cm} Q_1, Q_2, Q_3, \hspace{0.2cm} Q_1, Q_2, Q_3,
Q_4, \hspace{0.2cm} Q_1, Q_2, \ldots \]

At each time $t$, let \red{$m_i(t)$} be the smallest $m$
\green{(if it exists)} such that \green{the output
\red{$y_{i,m}(t)$} of $Q_m$ at node $i$}  is a \green{singleton}
(not an interval). \green{We identify this singleton with the
numerical value of its single element, and we} set $y_i(t) = y_{i,
m_i(t)}(t)/m_i(t)$. If \green{$m_i(t)$ is undefined},  then
$y_i(t)$ is set to some default value, \jnt{e.g., $\emptyset$.}

\jnt{Let us fix a value of $n$. For any $m\leq n$,}
the definition of $Q_m$ and Theorem
\ref{thm:average_computation} \jnt{imply} that there exists a time after which
the outputs $\red{y_{i,m}}$ of \jnt{$Q_m$}
do not change, and are \jnt{equal to a common value, denoted $y_m$, for every $i$.} Moreover, at least one of \jnt{the
algorithms $Q_1,\ldots,Q_n$} has an integer output $y_m$. Indeed,  observe that $Q_n$
computes $(1/n) \sum_{i=1}^n n1_{x_i=1} = \sum_{i=1}^n 1_{x_i=1}$,
which is clearly an integer. \red{In particular, $m_i(t)$ is
eventually well-defined and bounded above by $n$.} We conclude
that there exists a time after which the output \jnt{$y_i(t)$} of our
\green{overall} algorithm is fixed, shared by all nodes, and
different from the default value \jnt{$\emptyset$}.

We now argue that this value is indeed $p_1$. \red{Let
$m^*$ be} the smallest $m$ for which the eventual output of $Q_m$
is a single integer $y_m$. Note that $y_{m^*}$  is the exact
average
of the $x_{i,m^*}$, i.e., \[ y_{m^*} = \frac{1}{n} \sum_{i=1}^n
m^* 1_{x_i=1} = m^* p_1.\] \red{For large $t$, we have
\jnt{$m_i(t)=m^*$ and therefore} $y_i(t) =
y_{i, m^*}(t)/m^* =p_1$, \jnt{as desired.}}

Finally, it remains to argue that the algorithm described here can
be implemented with a sequence of \jnt{infinite memory} automata. All the above
algorithm does is run a copy of all the automata implementing
$Q_1, Q_2, \ldots$ with time-dependent transitions.  This can be
accomplished with an automaton whose state space is \green{the
countable set} $\mathcal{N} \times \red{\cup_{m
=1}^{\infty}\prod_{i=1}^{m} {\cal Q}_i}$, where ${\cal Q}_i$ is
the state \green{space} of $Q_i$, and the set $\mathcal{N}$ of
integers is used to keep track of time.
\end{proof}

%% file: maxtracking.tex
\section{Computing and tracking maximal values}\label{sec:max_tracking}
\rt{We now describe an algorithm that tracks the maximum (over all
nodes) of time-varying inputs at each node. It will be used as a
subroutine of the \jnt{interval-averaging} algorithm described in
Section \ref{sec:comput_average}, \jnt{and \aoc{which is} used} to prove Theorem
\ref{thm:average_computation}.} \jnt{The basic idea is the same as
for  the simple algorithm for the detection problem given in Section
\ref{sec:formal_descr}: every node keeps track of the largest value
it has heard so far, and forwards this ``intermediate result'' to
its neighbors. However, when an input value changes, the existing
intermediate results \aoc{need to be} invalidated, and this is done by sending
``restart'' messages. A complication arises because invalidated
intermediate results might keep circulating in the network, always one
step ahead of the restart messages. We deal with this difficulty by
``slowing down'' the intermediate results, so that they travel at
half the speed of the restart messages.} \jblue{In this manner,
restart messages are guaranteed to eventually catch up with and
remove invalidated intermediate results.}

\jnt{We start by giving the specifications of the algorithm.}
Suppose that each node $i$ has \green{a time-varying input} $u_i(t)$
stored in memory at time $t$, \green{belonging to a finite} set of
numbers ${\cal U}$. \green{We assume that, for each $i$, the
sequence $u_i(t)$} must eventually stop changing, i.e., \green{that
there exists} some $T'$ \green{such that}
\[ u_i(t) = u_i(T'), ~~~~~~~~\mbox{ for all  } \jnt{i \mbox{ and }}t \geq T'.\]
(However, \jnt{node $i$} need not \jnt{ever be} aware that $u_i(t)$
has reached its final value.) Our goal is to develop a distributed
algorithm \green{whose output eventually settles on the value
$\max_i u_i(T')$. More precisely, each} node $i$ is to maintain a
number $M_i(t)$ which must satisfy the following \jblue{condition:}
\red{for every \jnt{network}} \green{and any allowed sequences
$u_i(t)$,} there exists some $T''$ with
\[ M_i(t) = \max_{\jnt{j}=1,\ldots,n} u_j(t), ~~~~~\mbox{  for all  }
\jnt{i \mbox{ and }}t \geq T''. \]

Moreover, \jnt{each} node $i$ must also maintain a pointer $P_i(t)$
to a neighbor or to itself. We will use the notation
$P^2_i(t)=P_{P_i(t)}(t)$, $P^3_i(t)=P_{P_i^2(t)}(t)$, etc. \green{We
require} the following \green{additional} property, for all $t$
larger than $T''$: for each node $i$ there exists a node $\green{j}$
and a power $K$ such that for all $k \geq K$ we have $P^k_i\aoc{(t)} =  j$ and $M_i(t) = u_j(t)$. In  words, by successively following
the pointers $P_i(t)$, one \green{can} arrive at a node with a
maximal value.

 \jmj{We} next describe the algorithm. \jnt{We will use the term {\em slot $t$} to refer, loosely
speaking, to the interval between times $t$ and $t+1$. More
precisely, during slot $t$} each node processes the messages that
have arrived at time $t$ and computes the state at time $t+1$ as
well as the messages it will send at time $t+1$.

\begin{figure}[h] \hspace{0cm}\centering
 \epsfig{file=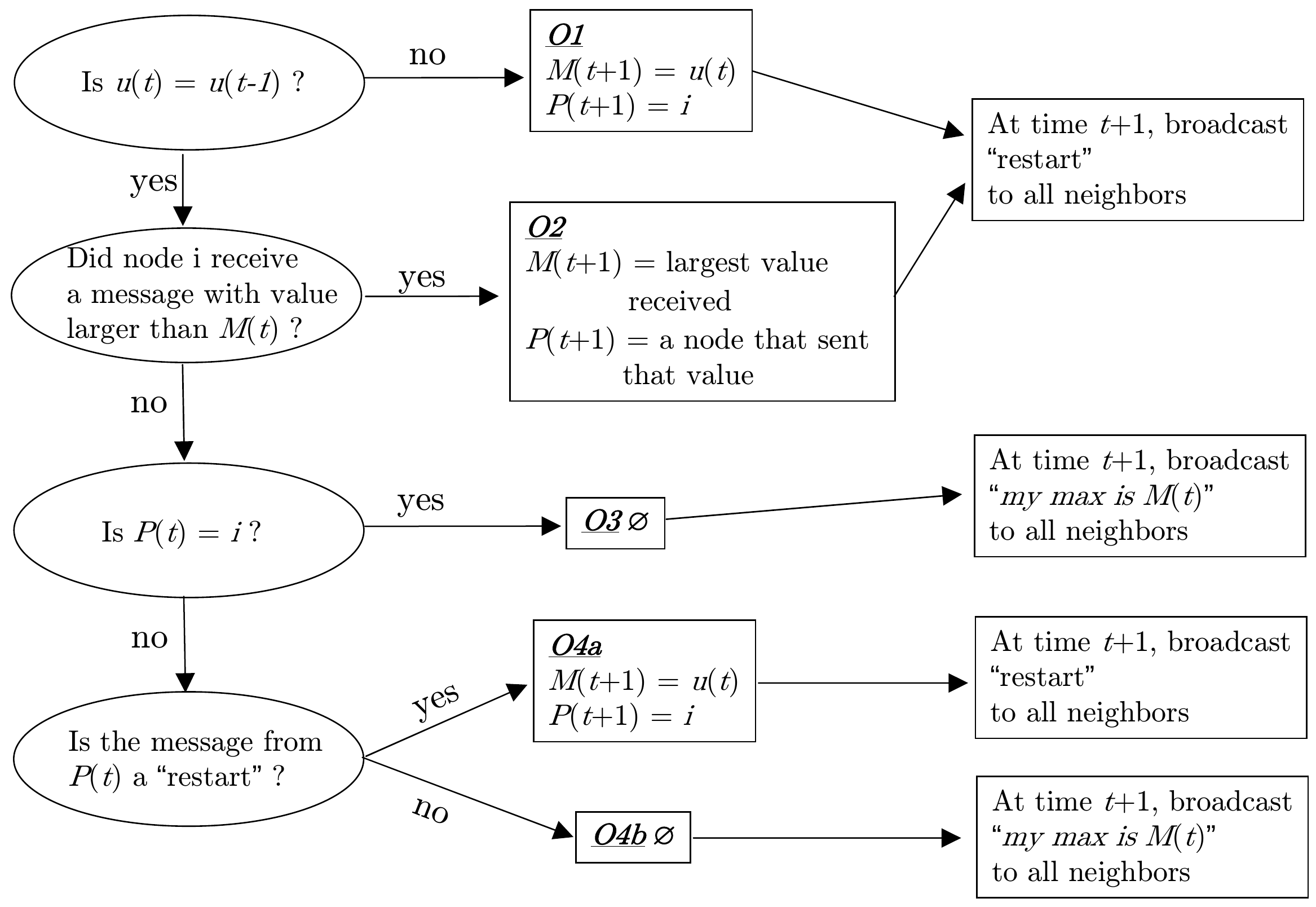, scale=.6}
\caption{\jnt{Flowchart of the procedure used by node $i$ during
slot $t$ in the maximum tracking algorithm. The subscript $i$ is
omitted, but $u(t)$, $M(t)$, and $P(t)$ should be understood as
$u_i(t)$, $M_i(t)$, and $P_i(t)$. In those cases where an updated
value of $M$ or $P$ is not indicated, it is assumed that
$M(t+1)=M(t)$ and $P(t+1)=P(t)$. \jblue{The symbol $\emptyset$ is
used to indicate no action. \jred{Note that the various actions
indicated are taken during slot $t$, but the messages determined by
these actions are sent (and instantaneously received) at time
$t+1$.}
 }}\jmj{Finally, observe that every node sends an \aoc{identical}  message to all its
neighbors at every time $t>0$.} \aoc{We note that the apparent non-determinism in instruction
O2 can be removed by picking a node with, say, the smallest port label.} }\label{fig:maxtracking}
\end{figure}

The variables $M_i(t)$ and $P_i(t)$ are a complete description of
the state of node $i$ at time $t$. Our algorithm has only two types
of messages \jnt{that a node} can send to \jnt{its} neighbors. Both
are broadcasts, in the sense that the node sends them to every
neighbor:
\begin{enumerate}
\item[1.]
``Restart!'' \item[2.] ``My \jnt{estimate of the maximum} is $y$,''
where $y$ is some number in \jnt{${\cal U}$} \jnt{chosen} by the
node.
\end{enumerate}
Initially, each node sets $M_i(0)=\rt{u}_i(0)$ and $P_i(0)=i$. At
time $t=0,1,\ldots$, nodes exchange messages, which \jnt{then}
determine their state at time $t+1$, i.e., the pair  $M_i(t+1),
P_i(t+1)$, \jnt{as well as the messages to be sent at time $t+1$.}
The procedure node $i$ uses to do this is \jnt{described} in
Figure~\ref{fig:maxtracking}. One can verify that \jnt{a memory size
of $C\log |{\cal U}|+C\log d(i)$ at each node $i$ suffices,
\jblue{where $C$ is an absolute constant.} (This is because $M_i$
and $P_i$ can take one of $|{\cal U}|$ and $d(i)$ possible values,
respectively.)}

\aoe{The result that follows asserts the correctness of the
algorithm. The idea of the proof is quite simple. Nodes maintain
estimates $M_i(t)$ which track the largest among all the $u_i(t)$ in
the graph; these estimates are ``slowly" forwarded by the nodes to
their neighbors, with many artificial delays along the way. Should
some value $u_j(t)$ change, restart messages traveling without
artificial delays  are forwarded to any node which thought $j$ had a
maximal value, causing those nodes to start over. The possibility of
cycling between restarts and forwards is avoided because restarts
travel faster. Eventually, the variables $u_i(t)$ stop changing,
and the algorithm settles on the correct answer.} \aoe{On the other hand a formal \aoc{and
rigorous exposition} of this simple argument is
rather tedious. \aoc{A formal proof is available in 
Appendix A.
}}

\begin{theorem}[Correctness of \jnt{the} maximum tracking algorithm]\label{thm:maxtracking}
\jnt{Suppose that the $u_i(t)$ stop changing after some finite time.
Then, for every network,} there is a time after which \jnt{the
variables $P_i(t)$ and $M_i(t)$ stop changing and satisfy
$M_i(t) = \max_i \rt{u}_i(t)$; furthermore, after that time, and for
every $i$, the node $j=P_i^n(t)$ satisfies
$M_i(t)=u_j(t)$.}\end{theorem}

%% file: average_computation.tex
\section{ \jnt{Interval-averaging}}\label{sec:comput_average}

\jnt{In this section, we present an interval-averaging algorithm
and prove its correctness. We start with an informal discussion of
the main idea.} Imagine the integer \jnt{input} value $x_i$
\green{as represented by a number of $x_i$} pebbles \jblue{at node
$i$.} \jnt{The} algorithm attempts to exchange pebbles between
nodes with unequal numbers so that the overall distribution
\green{becomes} more even. Eventually, either all nodes will have
the same number of pebbles, or some will have a certain number and
others just one more. \red{We let $u_i(t)$ be the current number
of pebbles at node $i$; in particular, $u_i(0)=x_i$. An important
property of the algorithm \ora{will be} that the total number of
pebbles is conserved.}

To match nodes with unequal number of pebbles, we use the maximum
tracking algorithm of Section \ref{sec:max_tracking}. Recall that
the algorithm provides nodes with pointers which attempt to track
the \green{location of the} maximal values. When a node with
$\red{u}_i$ pebbles comes to believe in this way that a node with
at least $u_i+2$ pebbles exists, it sends a request in the
direction \green{of the latter node} to obtain one or more
pebbles. This request follows \blue{a} path to \blue{a} node with
\blue{a} maximal number of pebbles until \blue{the request}
either gets denied, or gets accepted by a node with at least
$u_i+2$ pebbles.

\subsection{\jnt{The algorithm} \label{se:alg}}
\jnt{The} algorithm uses two types of messages. \jblue{Each type
of message}  \jnt{can be either \emph{originated} at a node or
\emph{forwarded} by a node.}

\smallskip
\texitem{(a)} \green{(Request, $r$): This is a request for a
transfer of \jnt{pebbles.} Here, $r$ is an integer that} represents the
number of pebbles \ora{$u_i(t)$} at the
\jnt{node \ora{$i$} that first originated the request,}
\jblue{at the time \ora{$t$} that the request was originated.}
\ora{(Note, however, that this request is actually sent at time $t+1$.)}
\smallskip
\texitem{(b)}
\green{(Accept, $w$): This corresponds to} acceptance of
a request, and \jblue{a} transfer of  $w$ pebbles
\jnt{towards the node that originated the request.}
\jnt{An acceptance} with a value $w=0$ represents a request denial.
\smallskip

\blue{As} part of the algorithm, the nodes run the maximum
tracking algorithm of Section \ref{sec:max_tracking}, as well as a
minimum tracking counterpart. In particular, each node $i$
has access to the variables $M_i(t)$ and $P_i(t)$ of the maximum
tracking algorithm (recall that these are, respectively, the
\green{estimated maximum and a} pointer to a neighbor \aoc{or to itself}).
\green{Furthermore,} each node maintains three additional
variables.

\smallskip
\texitem{(a)}
``\jnt{M}ode($t$)"$\in$ \{\jnt{F}ree,\ \jnt{B}locked\}. Initially, the mode of every
node is free.
\jnt{A node is blocked if it has originated or forwarded a request,
\jblue{and is still waiting to hear whether the request is} accepted (or denied).}

\smallskip
\texitem{(b)} ``$\Rin_i(t)$'' and ``$\Rout_i(t)$" are pointers to
a neighbor of $i$, or to \jnt{$i$} itself. \jnt{The \ora{meaning}
of these pointers \ora{when in blocked mode} are as follows. If
$\Rout_i(t)=j$, then node $i$ has sent (either originated or
forwarded) a request to node $j$, \jblue{and is still in blocked
mode, waiting to hear whether the request} \ora{is acceptd or
denied.} If $\Rin_i(t)=k$, and $k\neq i$, then node $i$ has
received a request from node $k$ \jblue{but has not yet responded
to node $k$.} If $\Rin_i(t)=i$, then node $i$ has originated a
request \jblue{and is still in blocked mode, waiting to hear
whether the request} is accepted or denied.}

\smallskip

\begin{figure}
\centering

\includegraphics[scale = .50]{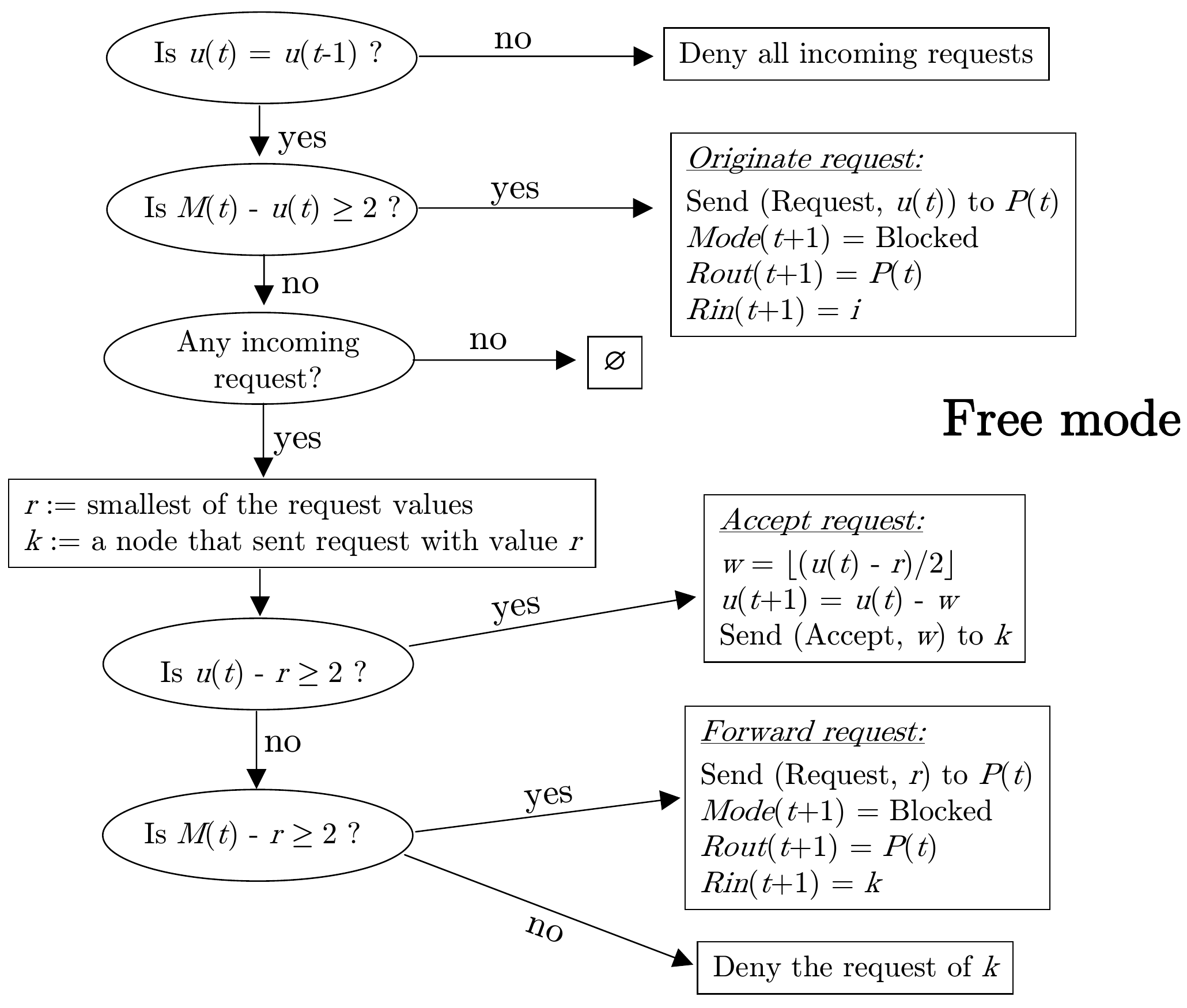}\\ 

\vspace{30pt}

\includegraphics[scale = .55]{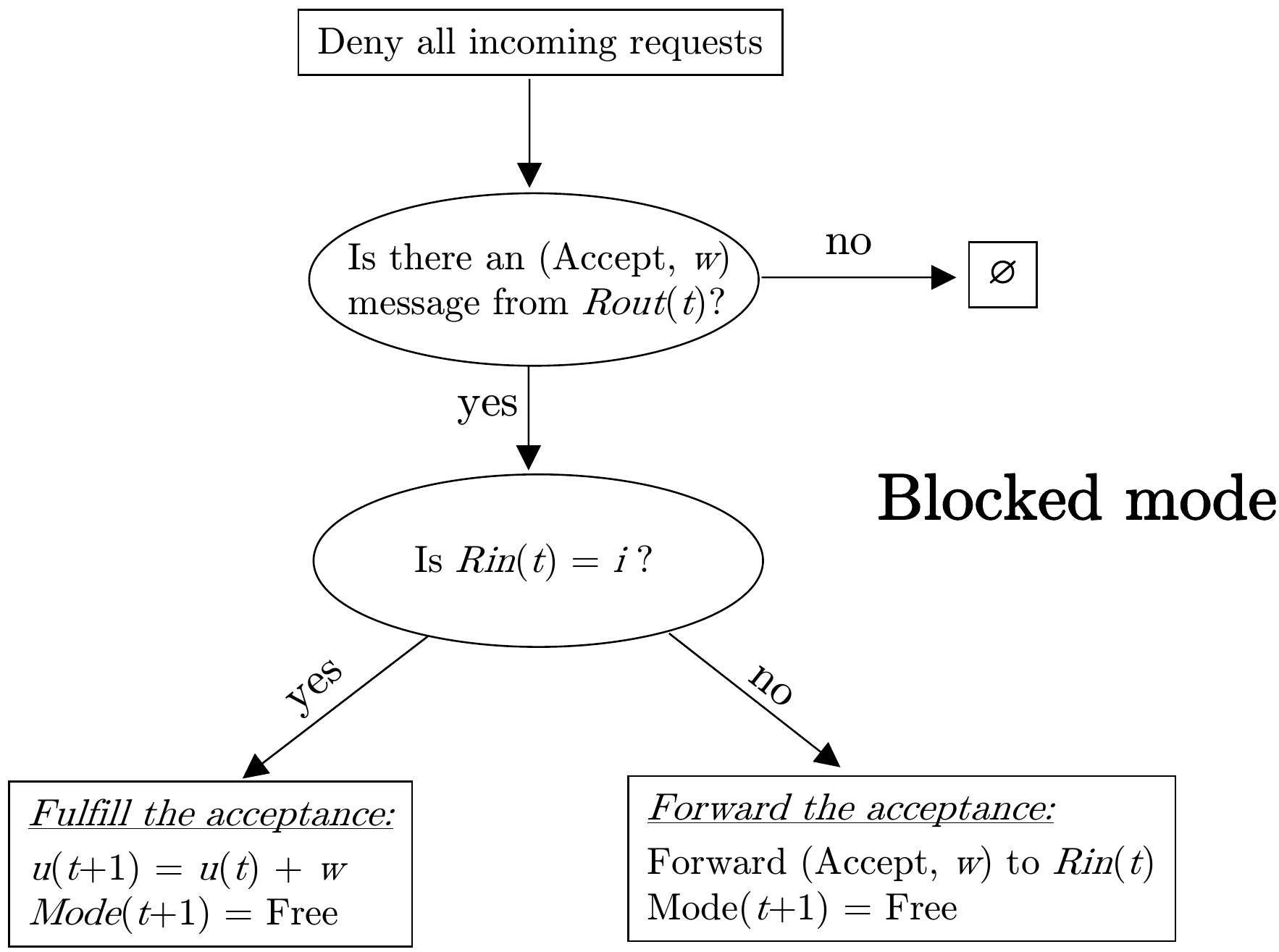}  

\caption{\jnt{Flowchart of the procedure used by node $i$ during slot $t$ in the interval-averaging algorithm.
The subscript $i$ is omitted from variables such as Mode$(t)$, $M(t)$, etc.
Variables for which an update is not explicitly indicated are assumed to remain unchanged.}
\jblue{``Denying a request'' is a shorthand for
$i$ sending a message of the form (Accept, 0) at time $t+1$ to a node from which $i$ received a request at time $t$. Note also that \ora{``forward the acceptance''} in the blocked mode includes the case where
the answer had $w=0$ (i.e., \ora{it was} a request denial), in which case the denial is forwarded.}
}\label{fig:descr_avg_algo}
\end{figure}

\jnt{A precise description of the algorithm is given in Figure
\ref{fig:descr_avg_algo}. The proof of  correctness is given in the
next subsection, thus also establishing Theorem
\ref{thm:average_computation}. Furthermore,  we will show that the
time until the algorithm settles on the correct output is of order
$O(n^2K^2)$.}

\subsection{\jnt{Proof of correctness}}

\old{\jnt{To facilitate understanding, let us start with some observations.
A node alternates between the Free and Blocked modes. Whenever
it switches from Blocked to Free, $\Rin_i$ and $\Rout_i$ are set
to $\emptyset$; and whenever it switches from Free to Blocked,
$\Rin_i$ and $\Rout_i$ are set to values other than $\emptyset$.
Thus, Mode$_i(t)=$Free if and only if $\Rout_i(t)=\Rin_i(t)=\emptyset$.}}

We begin by arguing that that the rules of the algorithm
preclude one potential obstacle; \ora{we will show that} nodes will not get stuck sending
requests to themselves.

\begin{lemma}\label{l:no-self} \jnt{A node never sends (originates or forwards) a request to itself. More precisely, $\Rout_i(t)\neq i$, for all $i$ and $t$.}
\end{lemma}
\begin{proof} \jnt{By inspecting} the first two \jnt{cases for the} free mode, we observe that \jnt{if} node $i$
\jnt{originates} a request \ora{during time slot $t$ (and sends a
request message at time $t+1$),} \aoe{ then $P_i(t) \neq i$. Indeed,
to send a message, it must be true $M_i(t) > u_i(t) = u_i(t-1)$.
However, any action of the maximum tracking algorithm that sets
$P_i(t)=i$ also sets $M_i(t)=u_i(t-1)$, and moreover, as long as
$P_i$ doesn't change neither does $M_i$.  So }
\jnt{the recipient $P_i(t)$ of the request originated by $i$ is
different than $i$, \ora{and accordingly,} $\Rout_i\jblue{(t+1)}$ is
set to a value different than $i$. \ora{We argue} that the same is
true for the case where $\Rout_i$ is set by the the ``Forward
request'' box of the free mode. Indeed, that box is enabled  only
when $u_i(t)=u_i(t-1)$ and $u_i(t)-1\leq r <M_i(t)-1$, so that
$u_i(t-1)<M_i(t)$. \ora{As in the previous case, this implies that}
$P_i(t)\neq i$ and that $\Rout_i(t+1)$ is again set to a value other
than $i$. We conclude that $\Rout_i(t)\neq i$ for all $i$ and $t$.}
\end{proof}

\ora{We will now analyze the evolution of the requests.
A request is originated at some time $\tau$ by some originator node $\ell$ who sets $\Rin_{\ell}(\tau+1)=\ell$ and sends the request to some node $i=\Rout_{\ell}(\tau+1)=P_{\ell}(\tau)$.
The recipient $i$ of the request either accepts/denies it, in which case $\Rin_i$ remains unchanged, or forwards it while also setting $\Rin_i(\tau+2)$ to $\ell$. The process then continues similarly. The end result is that at any given time $t$, a request initiated by node $\ell$ has resulted in a ``request path of node $\ell$ at time $t$,'' which
is a maximal sequence of nodes $\ell,i_1,\ldots,i_k$ with
$\Rin_{\ell}(t)=\ell$, $\Rin_{i_1}(t)=\ell$, and $\Rin_{i_m}(t)=i_{m-1}$ for
$m\leq k$.

\begin{lemma}\label{lem:descr_Gr}
At any given time, different request paths cannot intersect (they involve disjoint sets of nodes). Furthermore, at any given time, a request path cannot visit the same node more than once.
\end{lemma}
\begin{proof}
For any time $t$, we form a graph that consists of all edges that lie on some request path.
Once a node $i$ is added to some request path, and as long as that request path includes $i$, node $i$ remains in blocked mode and the value of $\Rin_i$ cannot change. This means that adding a new edge  that points into $i$ is impossible. This readily implies that cycles cannot be formed and also that two request paths cannot involve a
common node.
\end{proof}
}

\ora{We use $p_\ell(t)$ to denote the request path of node $\ell$ at time $t$, and
$s_\ell(t)$ to denote the last node
on this path.}
We will say that a request \jnt{originated} by node $\ell$ \emph{terminates} when
node $\ell$ receives an \jnt{(Accept, $w$)} message, with \jnt{any} value
$w$.

\begin{lemma}\label{lem:request_terminate}
Every request eventually terminates. Specifically, if node
\ora{$\ell$} \jnt{originates} a request at time $t'$ \ora{(and sends a request message at time $t'+1$),} then \jnt{there exists a}  later time
$t''\mod{\leq t'+n}$ \jnt{at which} node $s_r(t'')$ receives an ``accept
request'' message \jnt{(perhaps with $w=0$),} which is forwarded until it reaches $\ell$,
\mod{no later than \ora{time} $t''+n$}. 
\end{lemma}
\begin{proof} By the rules of our algorithm, node $\ell$ sends a request
message to node $P_{\ell}(t')$ at time \ora{$t'+1$.} If node $P_\ell(t')$ replies
at time \ora{$t'+2$}
with a ``deny request'' \ora{response}  to $\ell$'s request, then
\ora{the claim is true; otherwise,} observe that $p_\ell(t'+\ora{2})$ is nonempty and until
$s_\ell(t)$ receives an ``accept request'' message, the length of
$p_\ell(t)$ increases \jnt{at each time step.} Since
this length cannot be larger than $n\ora{-1}$, by Lemma \ref{lem:descr_Gr},
it follows \jnt{that} $s_\ell(t)$ receives an ``accept request'' message at most
$n$ steps after $\ell$ initiated the \ora{request.}
\jnt{One can then} easily show that this \jnt{acceptance} message is forwarded \jnt{backwards along the path} \ora{(and the request path keeps shrinking)}
until the acceptance message reaches $\ell$, at most $n$ steps later.
\end{proof}

\ora{The arguments so far had mostly to do with deadlock avoidance.
The next  lemma concerns the progress made by the algorithm. Recall
that a central idea of the algorithm is to conserve the total number
of ``pebbles,'' but this must include both pebbles possessed by
nodes and pebbles in transit. We capture the idea of ``pebbles in
transit'' by defining a new variable. If $i$ is the originator  of
some request path that is present at time $t$, and if the final node
$s_i(t)$ of that path} receives an \jnt{(Accept, $w$)} message at
time $t$, we let $w_i(t)$ be \jnt{the value $w$ in that message.
(This convention includes the special case where $w=0$,
corresponding to a denial of the request).} \ora{In all other cases,
we} set $w_i(t) =0$. Intuitively, $w_i(t)$ is the value that has
already been given away by a node who answered a request originated
by node $i$, and that will eventually be added to \ora{$u_i$}, once
the answer reaches \ora{$i$}.

We now define
$$\hat u_i(t)= u_i(t) + w_i(t).$$ By the rules of our algorithm,
if $w_i(t)\ora{=w}>0$, \ora{an amount $w$} will eventually be added to $u_i$, once
the \ora{acceptance message is forwarded back to $i$.}  The value $\hat u_i$ can thus be
seen as a future value of $u_i$, that includes its present value
and the \ora{value that has} been sent to $i$ but has not yet reached it.

\ora{The rules of our algorithm imply
that the sum of the $\hat u_i$ remain constant. Let $\bar x$ be the average of the initial values $x_i$. Then,
\[ \frac{1}{n} \sum_{i=1}^n \hat{u}_i(t) =
\frac{1}{n} \sum_{i=1}^n x_i =\bar{x}. \] We define the \emph{variance} function $V$ as
\[ V(t) = \sum_{i=1}^n (\hat{u}_i(t) - \bar{x})^2.\]}

\begin{lemma}\label{prop:finite_request_acceptances}
The number of times \ora{that} a node \jnt{can send an \ora{acceptance} message (Accept, $w$) with $w\neq 0$,}
is finite.
\end{lemma}
\begin{proof} Let us first describe the idea behind the proof. Suppose that nodes
could instantaneously transfer value to each other. \ora{It is
easily checked that if} a node $i$ \ora{transfers} an amount $w^*$
to a node $j$ with $u_i-u_j\geq 2$ and $1\leq w^*\leq
\frac{1}{2}(u_i-u_j)$, the  variance $\sum_i (u_i- \ora{\bar x})^2$
decreases by at least \aoe{2}. \ora{Thus, there can only be} a
finite number of \ora{such transfers.} In our model, \ora{the
situation is more complicated because transfers}  are not immediate
\ora{and involve a process of requests and acceptances. A key
element of the argument is to realize that the algorithm can be
interpreted as if it only involved instantaneous exchanges involving
disjoint pairs of nodes.}

\ora{Let us consider the difference  $V(t+1)-V(t)$ at some typical time
$t$.
Changes in $V$ are solely due to changes in the $\hat u_i$.
Note that if a node $i$ executes the ``\aoc{fulfill} the acceptance'' instruction at time $t$, node $i$ was the
originator of the request
and the request path has length zero, so that it is also the final node on the path, and $s_i(t)=i$. According to our definition, $w_i(t)$ is the value $w$ in the message received by node $s_i(t)=i$.
At the next time step, we have $w_i(t+1)=0$ but $u_i(t+1)=u_i(t)+w$.
Thus, $\hat u_i$ does not change, and the function $V$ is unaffected.}

\ora{By inspecting the algorithm, we see that  a nonzero  difference
$V(t+1)-V(t)$ is possible only if some node $i$ executes the ``accept request'' instruction at slot $t$, with
some particular value $w^*\aoc{\neq 0}$,
in which case $u_i(t+1)=u_i(t)-w^*$. For this to happen, node $i$ received a message (Request, $r$) at time $t$ from a node $k$ for which $\Rout_k(t)=i$, and with $u_i(t) - r\geq 2$.
That node $k$ was the last node, $s_{\ell}(t)$, on the request path of some originator node $\ell$. Node $k$ receives an (Accept, $w^*$)
message at time $t+1$ and, therefore, according to our definition, this sets $w_\ell(t+1)=w^*$.}

It follows from the rules of our algorithm that
$\ell$ \ora{had} originated a request with value $r=u_{\ell}(t')$ at
some previous time $t'$. Subsequently, node $\ell$ entered the blocked
mode, preventing any modification of $u_{\ell}$, so that
$r=u_{\ell}(t) = u_{\ell}(t+1)$. Moreover, observe that $w_{\ell}(t)$ was 0 because \ora{by time $t$,} no node
had answered \ora{$\ell$'s} request.
\ora{Furthermore,} $w_i(t+1) = w_i(t) =0$ because
having a positive $w_i$ requires \ora{$i$} to be in blocked mode, preventing
the execution of ``accept request". It follows that
$$\hat
u_i(t+1) = u_i(t+1) = u_i(t) - w^* = \hat u_i(t) - w^*,$$ and
$$\hat
u_{\ell}(t+1) = r + w^* = \hat u_{\ell}(t) + w^*.$$
\ora{Using the update equation $w^*=\lfloor(u_i(t)-r)/2 \rfloor$, and the fact $u_i(t) - r\geq 2$,
we obtain} $$1\leq w^* \leq \frac{1}{2}(u_i(t) -
r) = \frac{1}{2}(\hat u_i(t) - \hat u_{\ell}(t)).$$ Combining
with the previous equalities, we have
$$
\hat u_{\ell}(t) +1 \leq  \hat u_{\ell}(t+1) \leq \hat u_i(t+1)
\leq \hat u_i(t) -1.
$$
\ora{Assume for a moment that node $i$ was the only one that
executed the ``accept request'' instruction at time $t$. Then, all
of the variables $\hat u_j$, for $j\neq i,\ell$, remain unchanged.}
Simple algebraic manipulations then show that $V$ decreases by at
least \aoe{2}. \ora{If there was another pair of nodes, say $j$ and
$k$, that were involved in a transfer of value at time $t$, it is
not hard to see that the transfer of value was related to a
different request, involving a separate request path. In particular,
the pairs $\ell, i$ and $j,k$ do not overlap. This implies that the
cumulative effect of multiple transfers on the difference
$V(t+1)-V(t)$ is the sum of the effects of individual transfers.
Thus, at every time for which at least one ``accept request'' step
is executed, $V$ decreases by at least 2. We also see that no
operation can ever result in an increase of $V$. It follows that}
the instruction ``accept request" can be executed only a finite
number of times.
\end{proof}

\begin{proposition}\label{prop:final_proof_average} There is a time $t'$ such that $\rt{u}_i(t) =
\rt{u}_i(t')$,  for \jnt{all $i$ and} all $t \geq t'$. Moreover,
\begin{eqnarray*}
\sum_\jnt{i} \rt{u}_i(t') & = & \sum_\jnt{i} x_i, \\
\jnt{\max_i \rt{u}_i(t')-\min_i \rt{u}_i(t')} & \leq &   1.
\end{eqnarray*}
\end{proposition}

\begin{proof}
\noindent It follows from \jnt{Lemma}
\ref{prop:finite_request_acceptances} that there is a time \aoc{$t'$}after
which no more requests are accepted \jnt{with $w\neq 0$.} By Lemma \ref{lem:request_terminate}, this
implies that after at most $n$ additional time steps, the system
will never \ora{again} contain any  ``accept request'' messages \jnt{with $w\neq 0$,} so no node
will change its value \ora{$u_i(t)$ thereafter.}

We have already argued that the \ora{sum (and therefore the
average) of the variables $\hat u_i(t)$} does not change.
\jnt{Once there are} no more ``accept request'' messages in the
system \jnt{with $w\neq 0$, we must have  $w_i(t)=0$, for all $i$.
Thus,} at this stage the \jnt{average of the $\rt{u}_i(t)$} is the
same as the \jnt{average} of \jnt{the $x_i$.}

It remains to show that once \jnt{the $\rt{u}_i(t)$ stop} changing, the
maximum and minimum \jnt{$u_i(t)$} differ by \jnt{at most} $1$. \jnt{Recall (cf.\ Theorem that
\ref{thm:maxtracking})} \aoc{that} at some time after \jnt{the $\rt{u}_i(t)$
stop} changing, \ora{all
estimates $M_i(t)$ of the maximum will be} equal to $M(t)$, the true maximum of \jnt{the
$\rt{u}_i(t)$;} moreover, \jnt{starting at
any node and} following the pointers $P_i(t)$  \jnt{leads to a node $j$} whose value \jnt{$u_j(t)$} is the true \jnt{maximum,} $M(t)$.
Now let \ora{$A$} be the set of nodes whose value at this stage is at
most $\max_i \rt{u}_i(t)-2$. To \ora{derive} a contradiction, \jnt{let us} suppose
that \ora{$A$} is nonempty.

\jnt{Because} only nodes in \ora{$A$} will \jnt{originate} requests, and \jnt{because} every
request eventually terminates \jnt{(cf.\ Lemma \ref{lem:request_terminate})}, if we wait
some finite amount of time, we will have
the additional property that all requests in the system originated
from \ora{$A$}. Moreover,  nodes in \ora{$A$}  \jnt{originate} requests every time they \aoc{are in} \ora{the
free mode,}
which is infinitely
often.

\ora{Consider now a request originating at a node in the set
\ora{$A$}. The value $r$ of such a request satisfies $M(t)-r\geq 2$,
which implies that every node that receives \aoc{it} either accepts it
(contradicting the fact that no more requests are accepted after
time \aoc{$t'$}), or forwards it, or denies it. But a node $i$ will deny a
request only if it is in blocked mode, that is, if it has already
forwarded some other request to node $P_i(t)$. This shows that
requests will keep propagating along links of the form $(i,P_i(t))$,
and therefore will eventually reach a node at which $u_i(t)=M(t)\geq
r+2$, at which point they will be accepted---a contradiction.}
\end{proof}

We are now ready to conclude.

\begin{proof}[Proof of Theorem \ref{thm:average_computation}]
\jnt{Let $u_i^*$ be the value that $u_i(t)$ eventually settles on.
Proposition \ref{prop:final_proof_average} \ora{readily} implies} that if the average $\bar x$ of
the $x_i$ is an integer, then $\rt{u}_i(t)=\rt{u}_i^*=\bar x$ will
eventually hold for every $i$. If \jnt{$\bar x$} is not an integer, then some
nodes will eventually have $\rt{u}_i(t)=\rt{u}_i^* = \lfloor \bar x
\rfloor$ and \ora{some other nodes} $\rt{u}_i(t) = \rt{u}_i^*=\lceil \bar x
\rceil$. Besides, using the maximum and minimum computation
algorithm, nodes will eventually have a correct estimate of $\max
\rt{u}_i^*$ and $\min \rt{u}_i^*$, since all
$\rt{u}_i(t)$ \jnt{settle on the fixed values $\rt{u}_i^*$.} This allows \jnt{the nodes} to \jnt{determine whether}
the average is exactly $\rt{u}^*_i$ (integer average), or \jnt{whether}  it
lies in $(\rt{u}_i^*,\rt{u}_i^*+1)$ or $( \rt{u}_i^*-1,
\rt{u}^*_i)$ (fractional average). Thus, \jnt{with some simple post-processing at each node (which can be done using finite automata), the nodes can produce the correct output for the interval-averaging problem. The proof of
Theorem \ref{thm:average_computation} is complete.}
\end{proof}

\aoe{Next, we give a convergence time bound for the algorithms we
have just described.}

\begin{theorem} \label{averagingtime}
\aoc{Any function $h$ satisfying the assumptions of
Theorem \ref{bounded-alg} can be computed in $O(n^2 K^2)$
time steps.}
\end{theorem}

\begin{theorem} \label{approximatetime}
\aoc{The functions $h_{\epsilon}$ whose existence is guaranteed by
Corollary \ref{cor:approx_set_epsilon} or Corollary
\ref{cor:approx_continuous} can be computed in time which grows
quadratically in $n$.  } \end{theorem}

\smallskip

\aoe{The general idea \aoc{behind Theorems \ref{averagingtime} and \ref{approximatetime}
is} quite simple. We have just argued that the
nonnegative function $V(t)$ decreases by at least $2$ each time a
request is accepted. \aoc{It also satisfies $V(0)=O(n K^2)$}. Thus there are at
most \aoc{$O(n K^2)$} acceptances. To prove \aoc{Theorems \ref{averagingtime} and \ref{approximatetime}
is}, one needs to
argue that if the algorithm has not terminated, there will be an
acceptance within $O(n)$ time steps. This should be \aoc{fairly clear} from the
proof \aoc{of Theorem \ref{thm:average_computation}. A formal argument is
given in 
Appendix B.
It is also shown there that the running time of our algorithm, for  many graphs,
satisfies a $\Omega(n^2)$ lower bound, in the worst case  over all initial conditions.  } }

%% file: simulations.tex
\section{Simulations \label{se:simul}} \begin{figure}
\begin{minipage}[b]{0.5\linewidth} 
\centering
\includegraphics[width=7cm]{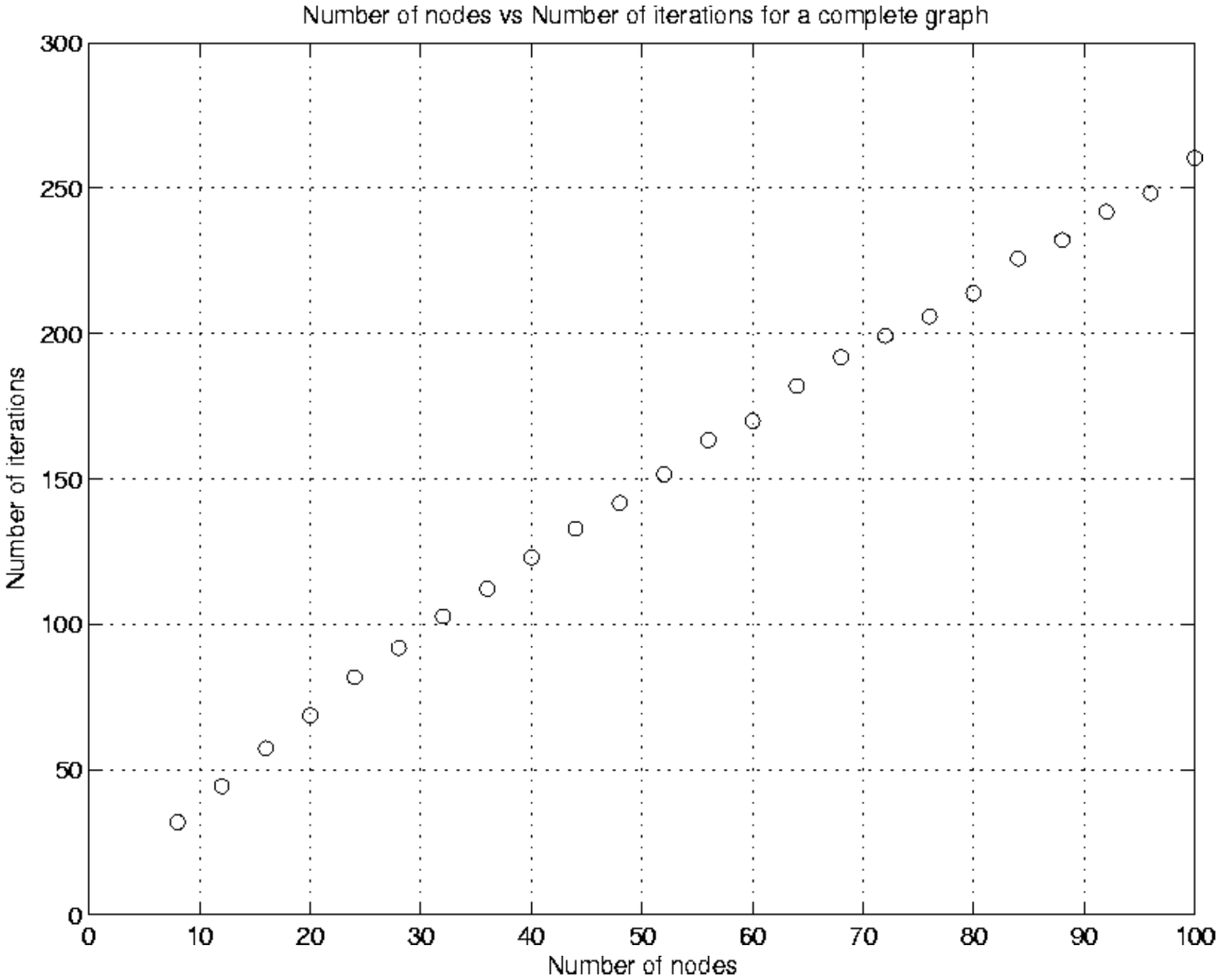}
\caption{\aoc{The number of iterations as a
function of the number of nodes for a complete graph under random initial conditions. \label{simulc}}}
\end{minipage}
\hspace{0.5cm} 
\begin{minipage}[b]{0.5\linewidth}
\centering
\includegraphics[width=7cm]{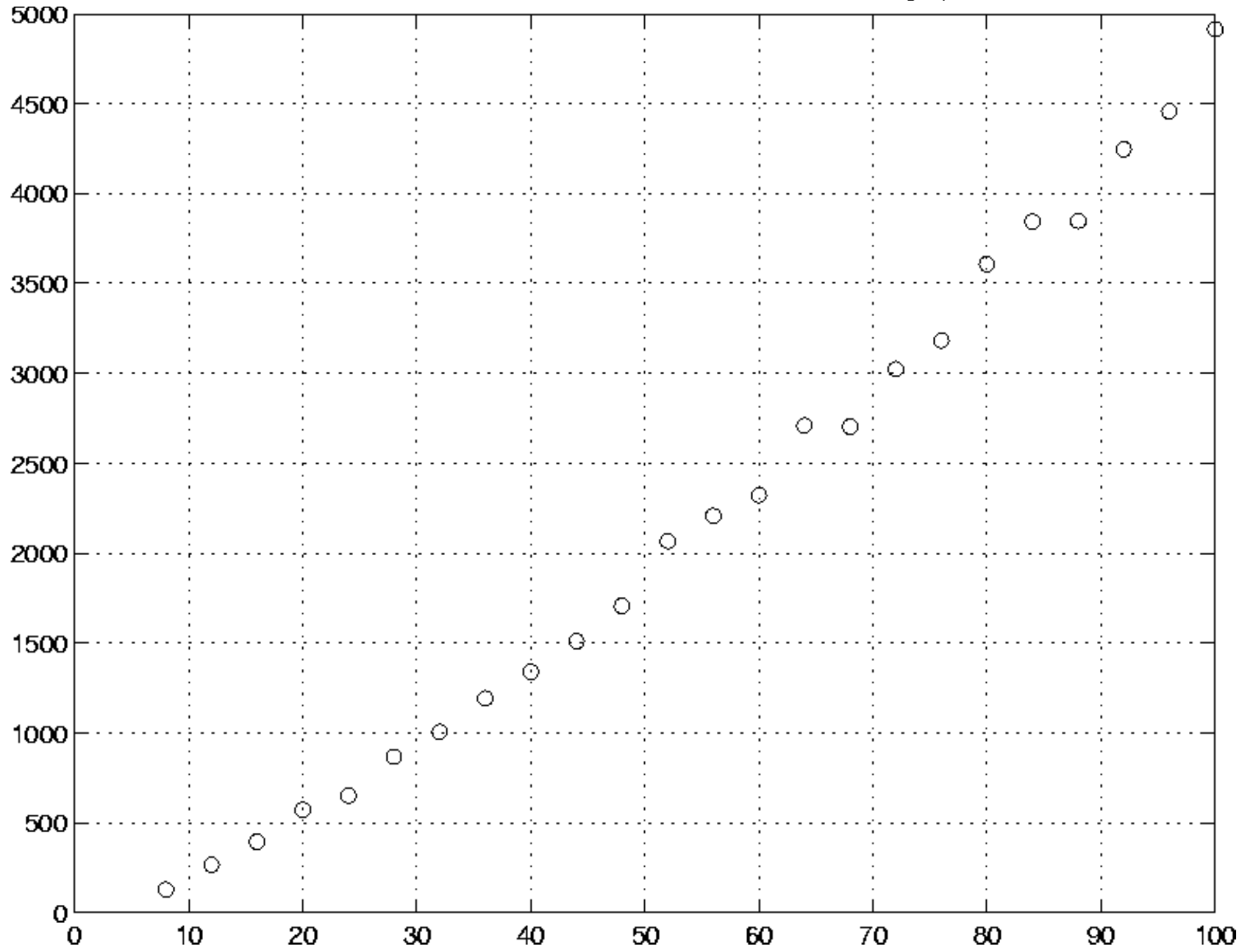}
\caption{\aoc{The number of iterations as a function
of the number of nodes for a line graph under random initial conditions. \label{simull}}}
\end{minipage}
\end{figure} \aoe{We report here on simulations involving our
algorithm on several natural graphs. \aoc{Figures \ref{simulc} and
\ref{simull} describe the results for a complete graph and a line. Initial conditions were
random integers between $1$ and $30$, and each data point represents
the average of two hundred runs.} As expected, convergence is faster on the complete graph. Moreover,
\aoc{convergence time in both simulations appears to be approximately linear.}

\aoe{Finally recall that it is possible for our algorithm to take on the order of $n^2$ (as  opposed to
$O(n)$) time \rev{steps} to converge. } \jntr{Figure \ref{dumbbellsimul} shows  simulation results} for the dumbbell graph (two complete graphs \jntr{with $n/3$ nodes,} 
connected by \aoc{a line}) \jntr{of length $n/3$; each} node in \aoc{one of the
complete graphs} starts
with $x_i(0)=1$, every node \aoc{in the other complete graph} starts with $x_i(0)=30$.
The time to converge in this case is quadratic in $n$.
\begin{figure}[h] \hspace{0cm}\centering
 \epsfig{file=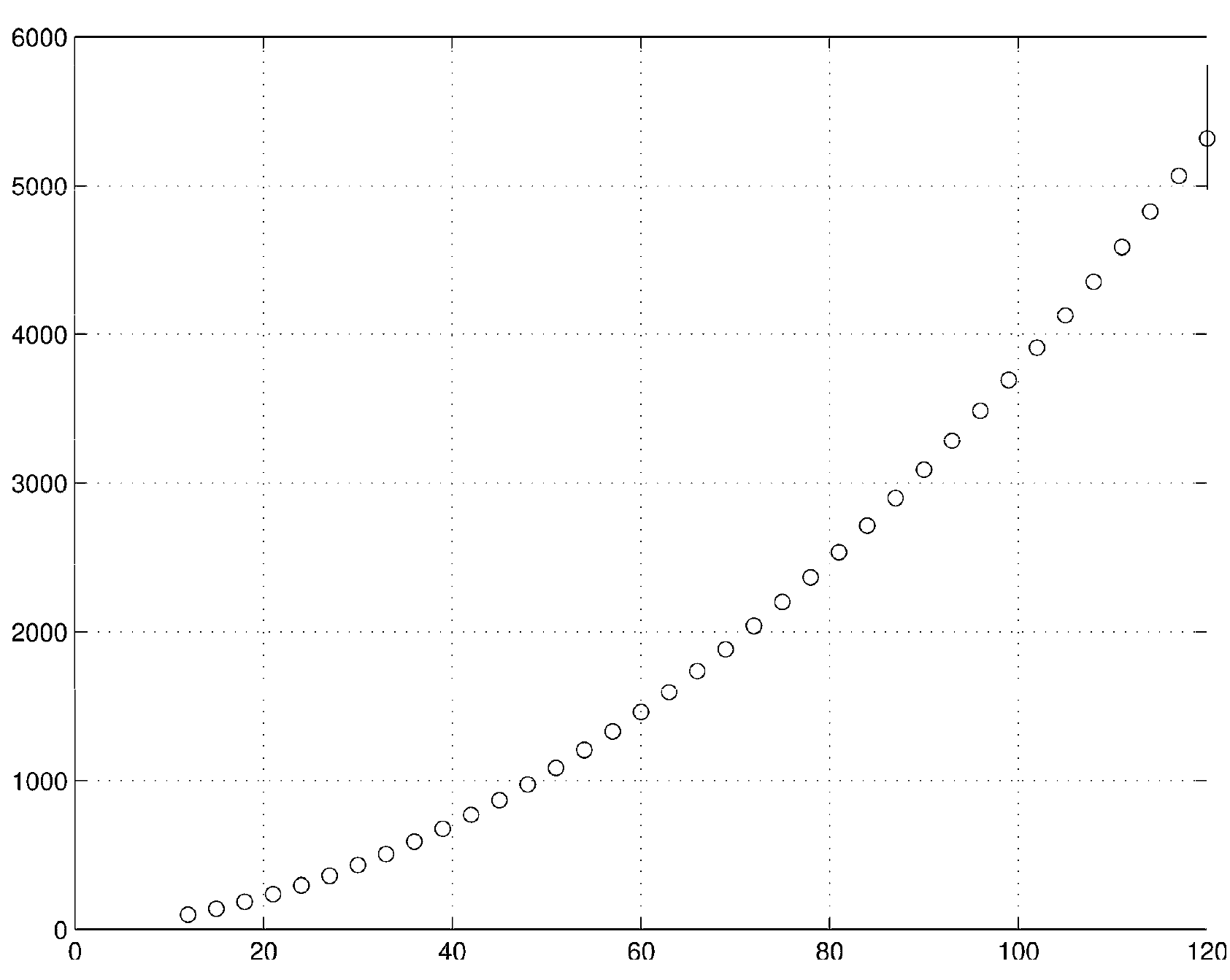, scale=.4} \caption{\aoc{The number of iterations as a function
of the number of nodes for a dumbbell graph.}} \label{dumbbellsimul}
\end{figure}}

%% file: conclusions.tex
\section{\rt{Conclusions}}\label{sec:conclusions}

We have proposed a model of \ora{deterministic anonymous distributed computation,} inspired by the \ora{wireless sensor network and} multi-agent
control literature. We have given an almost tight characterization
of the functions that are computable in our model. We have shown
that computable functions must depend only on the \ora{the frequencies with which}
the different initial conditions \ora{appear,} and that if this dependence can
be expressed in term of linear inequalities with rational
coefficients, the function is indeed computable. Under weaker
conditions, the function can be approximated with arbitrary
precision. It remains open to exactly characterize  the class of
computable function {\ora families}.

Our positive results are proved constructively, by providing \aoc{a}
generic algorithm for computing the desired functions. Interestingly,
the finite memory requirement is not used in our negative
results, which remain thus valid in the infinite memory case. In
particular, we have no examples of functions that can be computed
with infinite memory but are provably not computable with
finite memory. We suspect though that simple examples
exist; a good candidate could be the indicator function
$1_{p_1<1/\pi}$, which checks whether the fraction of nodes with a
particular initial condition is smaller than
$1/\pi$.

\ora{We have shown that our generic algorithms terminate in $O(n^2)$
time. On the other hand, it is clear that the termination time cannot
be faster than the graph diameter, which is of order $n$, in the worst case.}
Some problems, such as the detection problem described
in Section \ref{sec:formal_descr}, \ora{admit $O(n)$ algorithms.
On the other hand, it is an open problem whether
the interval averaging problem admits an $o(n^2)$ algorithm \aoc{under our model}.} \aoc{Finally,
we conjecture that the dependence on $K$ in our $O(n^2 K^2)$ complexity estimate can be
improved by designing a different algorithm.}

\ora{Possible extensions of this work involve variations of the model of computation.
For example, the algorithm for detection problem, described
in Section \ref{sec:formal_descr},
does not make full use of the flexibility allowed by our model of computation. In particular,
 for any given $i$, the messages $m_{ij}(t)$ are the same for all $j$, so,
in some sense, messages are ``broadcast'' as opposed to being personalized for
the different outgoing links.  This raises an interesting question:
 do there exist computable function families that
become non-computable when we restrict to algorithms that are
limited to broadcast messages? \aoc{We have reasons to believe
that in a pure broadcast scenario where nodes \aoc{located in a non-bidirectional
network} broadcast messages
without knowing their out-degree (i.e., the size of their
audience), the only computable functions are those which test
whether there exists a node whose initial condition belongs to
a given subset of $\{0, 1,\ldots,K\}$, and combinations of such
functions.}}

\ora{Another important direction is to consider models in which the underlying graph may vary \aoc{with} time.
It is of interest to develop algorithms that converge to the correct answer at least when the underlying graph eventually stops changing.
For the case where the graph keeps changing while maintaining some degree of connectivity, we conjecture that no deterministic algorithm with bounded memory can solve the interval-averaging problem.} \ora{Finally, other extensions involve models accounting for clock asynchronism,
delays in message propagation, \rev{or data loss.}
}
\vspace{-10pt}

%% file: biblio.tex
\begin{small}

\end{small}

%% file: appen_proof_max.tex
\def\aoe{}
\def\aon#1{{#1}}
\section{Proof of Theorem \ref{thm:maxtracking}}\label{appen:proof_max_track}

\aon{In the following, we will occasionally} use the following convenient shorthand:
we will say that a statement $S(t)$ {\em holds eventually} if
there exists some $T$ so that \jnt{$S(t)$ is} true for all $t \geq
T$.

The analysis is made easier by introducing the time-varying directed
graph $G(t) = (\{1,\ldots,n\}, E(t))$ where $(i,j) \in E(t)$ if
\jnt{$i\neq j$, } $P_i(t)=j$, and there is no ``restart" message
sent \jnt{by $j$ (and therefore received by $i$) during \aoe{time}
$t$.} \aon{We will abuse notation by writing $(i,j) \in G(t)$ to mean that
the edge $(i,j)$ belongs to the set of edges $E(t)$ at time $n$.}

\begin{lemma} \label{newedge} \jnt{Suppose that $(i,j) \notin G(t-1)$ and $(i,j) \in G(t)$.}
Then, $i$ executes O2 during slot $t-1$.
\end{lemma}

\begin{proof}
\jnt{Suppose that $(i,j) \notin G(t-1)$ and $(i,j) \in G(t)$.} If
$P_i(t-1)=j$, \jnt{the definition of $G(t-1)$ implies}  that $j$
sent a restart message to $i$ at time $t-1$. \jmj{Moreover, $i$
cannot execute O3 during the time slot $t-1$ as \aon{this} would
require $P_i(t-1) = i$.}\jmj{Therefore, during this time slot,}
\jnt{node $i$ either executes O2 (and we are
done\footnote{\jred{In fact, it can be shown that this case never
occurs.}}) or it executes one of O1 and O4a. For both of the
latter two cases, we will have} $P_i(t)=i$, so that $(i,j)$ will
not be in $G(t)$, \jnt{a contradiction.} Thus, it must be that
$P_i(t-1) \neq j$. We now observe that the only place in Figure
\ref{fig:maxtracking} \jnt{that can change $P_i(t-1) \aon{\neq j}$
to} $P_i(t)=j$ is O2.
\end{proof}

\begin{lemma}\label{lem:speed1/2}
In \jnt{each of the following three} cases, node $i$ has no incoming
edges in \aon{either graph $G(t)$ or $G(t+1)$}:\\
(a) Node $i$ executes O1, O2, or O4a during time slot $t-1$;\\
(b) $M_i(t) \not = M_i(t-1)$;\\
(c) For some $j$, $(i,j) \in G(t)$ but $(i,j) \notin G(t-1).$
\end{lemma}
\begin{proof} \ao{(a) If $i$ executes O1, O2, or O4a during slot $t-1$, then it sends
a restart \jblue{message} to each of its neighbors at time $t$.
Then, \jblue{for any neighbor $j$} of $i$, \jnt{the definition of
$G(t)$ implies that}  $(j,i)$ is not in  $G(t)$. Moreover, by Lemma
\ref{newedge}, in order for $(j,i)$ to be in $G(t+1)$, \jblue{node
$j$ must execute} $O2$ during slot $t$. But the execution of O2
during slot $t$ cannot result in the addition of the edge $(j,i)$
\jblue{at time $t+1$,} \jnt{because} \jmj{the message broadcast by
$i$ \aon{at time $t$} to its neighbors was a restart.} So, $(j,i) \notin G(t+1)$. }

\ao{(b)} \jnt{If} $M_i(t) \not = M_i(t-1)$, \jnt{then} $i$ executes
O1, O2, or O4a during slot $t-1$, so the \jnt{claim} follows from
part (a).

\ao{(c)} By Lemma \ref{newedge}, \jnt{it must be the case that
\jblue{node} $i$ executes} O2 during slot $t-1$, and part (a)
implies the \jnt{result}.
\end{proof}

\begin{lemma}\label{lem:acyclic}
The graph $G(t)$ \ao{is acyclic}, for \jnt{all} $t$.
\end{lemma}
\begin{proof}
The initial graph $G(0)$ does not contain a cycle. \ao{ Let $t$ be
the first time a cycle \jblue{is present}, and let $(i,j)$ be an
edge \jnt{in} \jblue{a} cycle that \jnt{is added} at time $t$,
\jnt{i.e., $(i,j)$ belongs to $G(t)$ but not $G(t-1)$.} Lemma
\ref{lem:speed1/2}\jmj{(c)} implies that $i$ has no incoming edges
\jnt{in $G(t)$,} so $(i,j)$ cannot be an edge of the cycle---a
contradiction. }
\end{proof}

\jblue{Note that every node has out-degree at most one, because
$P_i(t)$ is a single-valued variable. Thus, the acyclic graph $G(t)$
must be a forest, specifically,} \jnt{a collection of disjoint
trees, with all arcs of a tree directed so that they point towards a
root of the tree (i.e., a node with zero out-degree).} The next
lemma establishes that $M_i$ is constant on any path of $G(t)$.

\begin{lemma}\label{lem:edge->mi=mj}
If \jnt{$(i,j)\in G(t)$,} then  $M_i(t) = M_j(t)$.
\end{lemma}

\begin{proof} \ao{Let $t'$ be \jnt{a time when $(i,j)$ is added to the graph, or more precisely,
a time such that $(i,j)\in G(t')$ but $(i,j)\notin G(t'-1)$.} First,
we argue that the statement we want to prove holds at time $t'$.
Indeed, Lemma \ref{newedge} implies that during slot $t'-1$, node
$i$ executed $O2$, so that $M_i(t') = M_j(t'-1)$. Moreover,
$M_j(t'-1)=M_j(t')$, \jnt{because otherwise case (b) in} Lemma
\ref{lem:speed1/2} would imply that} $j$ has no incoming edges at
time \jnt{$t'$, contradicting our assumption that $(i,j) \in
G(t')$.}

\ao{Next, we argue that \jnt{the property $M_i(t)=M_j(t)$ continues
to hold,} \jblue{starting from time $t'$ and for} as long as $(i,j)
\in G(t)$. Indeed, \jblue{as long as} $(i,j) \in G(t)$, then
$M_j(t)$ \jblue{remains} unchanged, \jnt{by case (b) of} Lemma
\ref{lem:speed1/2}. To argue that $M_i(t)$ \jblue{also remains}
unchanged, we simply observe that in Figure \ref{rings}, every box
which leads to a change in $M_i$ also \jnt{sets $P_i$ either to $i$
or to the sender of a message with value strictly larger than}
$M_i$; this latter \jblue{message} cannot come from $j$ because as
we just argued, increases in $M_j$ lead to \jnt{removal} of the edge
$(i,j)$ from $G(t)$. So, changes in $M_i$ are also accompanied by
\jnt{removal} of the edge $(i,j)$ from $G(t)$.}
\end{proof}

\jnt{For the purposes of the next lemma, we use the convention
$u_i(-1)=u_i(0)$.}

\begin{lemma}\label{lem:m_geq_y}
If $P_i(t) = i$, then $M_i(t) = \rt{u}_i(t-1)$; \jnt{if $P_i(t) \neq
i$, then} $M_i(t) > \rt{u}_i(t\ao{-}1)$.
\end{lemma}

\begin{proof}
We prove this result by induction. \jnt{Because of} the convention
$\rt{u}_i(-1) = \rt{u}_i(0)$, \jblue{and the initialization
$P_i(0)=i$, $M_i(0)=u_i(0)$,} the result trivially holds at time
$t=0$. Suppose now that the result holds at time $t$.
\ao{\jnt{During time slot $t$, we have three possibilities}
\jblue{for node $i$:}
\begin{itemize}
\item[(i)] Node $i$ executes O1 or O4a. In this case,
$M_i(t+1)=\rt{u}_i(t), P_i(t+1)=i$, so the result holds at time
$t+1$. \item[(ii)] Node $i$ executes O2. In this case $P_i(t) \neq
i$ and $M_i(t+1) > M_i(t) \geq \rt{u}_i(t-1)=\rt{u}_i(t)$. \jnt{The
first inequality follows from the condition for entering step O2.
The second follows from the induction hypothesis. The} last equality
follows because if $\rt{u}_i(t) \neq \rt{u}_i(t-1)$, node $i$ would
have executed O1 rather than O2. So, once again, the result holds at
time $t+1$. \item[(iii)] Node $i$ executes O3 or O4b.  The result
holds at time $t+1$ because neither $\rt{u}_i$ nor $M_i$ changes.
\end{itemize}}
\end{proof}

In the sequel, we will use $T'$ to refer to a time after which all
the $\rt{u}_i$ are constant. The following lemma shows that, after
$T'$, the largest estimate does not increase.

\begin{lemma}\label{lem:nondecreasing}\aoe{
Suppose that at some time $t' > T'$ we have $\widehat{M}> \max_i
M_i(t')$. Then \begin{equation} \label{eq:nondecreasing}
\widehat{M}> \max_i M_i(t) \end{equation} for all $t \geq t'$.}
\end{lemma}

\begin{proof}
{ \aoe We prove Eq. (\ref{eq:nondecreasing}) by induction. By
assumption it holds at time $t=t'$. Suppose now Eq.
(\ref{eq:nondecreasing}) holds at time $t$; we will show \aon{that} it holds at
time $t+1$.}

Consider a node $i$. If it executes O2 during the slot $t$, it sets
$M_i(t+1)$ to the value contained in a message sent at time $t$ by
some node $j$. It follows from the rules of our algorithm that the
value in this message is $M_j(t)$ and therefore, $M_i(t+1) = M_j(t)
< \widehat{M}$.

Any operation other than O2 that modifies $M_i$ sets $M_i(t+1) =
u_i(t)$, {\aoe and since $u_i(t)$ does not change after time $T'$,
we have $M_i(t+1) = u_i(t-1)$. By Lemma \ref{lem:m_geq_y},
 $M_i(t) \geq u_i(t-1)$, so that $M_i(t+1) \leq
M_i(t)$. We conclude that $M_i(t+1) < \widehat{M}$ holds for this
case as well.}
\end{proof}

We now introduce some terminology used to specify whether the
estimate $M_i(t)$ held by a node has been invalidated or not.
Formally, we say that node $i$ has a \emph{\jnt{valid} estimate}
\jnt{at time $t$} if by following the path in $G(t)$ \jnt{that
starts} at $i$, \jnt{we}  eventually arrive at a node $r$ with
$P_r(t) = r$ and \jblue{$M_i(t)=u_r(t-1)$.} In any other case, we
say that a node has an \emph{\jnt{invalid} estimate} at time $t$.

\smallskip

{\bf Remark:} Because of the acyclicity property, a path
in $G(t)$, starting from a node $i$, eventually leads to a node $r$
with out-degree 0; it follows from Lemma \ref{lem:edge->mi=mj} that
$M_r(t) = M_i(t)$. Moreover, Lemma \ref{lem:m_geq_y} implies that if
$P_r(t) = r$, then \jnt{$M_i(t) =M_r(t)=\rt{u}_r(\jblue{t-1})$,} so
that the estimate is \jnt{valid.} Therefore, if $i$ has an
\jnt{invalid} estimate, the corresponding node $r$ must have $P_r(t)
\not = r$; \jnt{on the other hand,} since $r$ has out-degree 0 in
$G(t)$, \jnt{the definition of $G(t)$} implies that there is a
``restart" message from $P_r(t)$ to $r$ sent at time $t$.

\smallskip

{\aoe The following lemma gives some conditions which allow us to
conclude that a given node has reached a \aon{final} state.}

\begin{lemma}\label{lem:valid+max-->final} {\aoe
Fix some $t'>T'$ and let $M^*$ be the largest estimate at time $t'$,
i.e., $M^*=\max_i M_i(t')$. If $M_i(t')
= M^*$, and this estimate is valid, then  for all $t \geq t'$: \\
\noindent (a) $M_i(t) = M^*$, $P_i(t) = P_i(t')$, and node $i$ has a valid estimate at time $t$. \\
(b) Node $i$ executes either O3 or O4b at time $t$.}
\end{lemma}
\begin{proof} {\aoe We will prove this result by induction on $t$.
Fix some node $i$. By assumption, part (a) holds at time $t=t'$.
To show part (b) at time $t=t'$, we first argue that $i$ does not
execute O2 during the time slot $t$. Indeed, this would require
$i$ to have received a message with an estimate strictly larger
than $M^*$, sent by some node $j$ who executed O3 or O4b during
the slot $t-1$. In either case, $M^* < M_j(t-1) = M_j(t)$,
contradicting the definition of $M^*$. \aon{Because of the
definition of $T'$, $u_i(t)=u_i(t-1)$ for $t>T'$, so that $i$ does
not execute O1.} This concludes the proof of the base case. }

{\aoe Next, we suppose that our result holds at time $t$, and we will
argue that it holds at time $t+1$. If } $P_i(t) = i$, then $i$
executes O3 {\aoe during slot $t$}, so that $M_i(t+1) = M_i(t)$ and
$P_i(t+1) = P_i(t)$, {\aoe completes the induction step for this case. }

{\aoe It remains to consider the case where $P_i(t) = j \neq i$. }
It follows from the definition of a valid estimate that $(i,j) \in
E(t)$. \aon{Using the definition of $E(t)$, we conclude that}
there is no restart message sent from $j$ to $i$ at time $t$.
{\aoe By the induction hypothesis,}
during the slot $t-1$, $j$ has thus executed O3 or O4b, so that
$M_j(t-1) = M_j(t)$; {\aoe in fact, Lemma \ref{lem:edge->mi=mj}
gives that $M_j(t)=M_i(t)=M^*$}. Thus during slot $t$, $i$ reads a
message from $j=P_i(t)$ with the estimate $M^*$, {\aoe and
executes O4b, consequently leaving its $M_i$ or $P_i$ unchanged.}

{\aoe We finally argue that node $i$'s estimate remains valid.
This is the case because since we can apply the arguments of the
previous paragraph to every node $j$ on the path from $i$ to a node
with out-degree $0$; we obtain that all of these nodes both (i)
keep $P_j(t+1)=P_j(t)$ and (ii) execute O3 or O4b, and consequently
do not send out any restart messages. }
\end{proof}

\smallskip

\aoc{Recall (see the comments following Lemma \ref{lem:acyclic})
that $G(t)$ consists of a  collection of disjoint in-trees (trees
in which all edges are oriented  towards a root node).
Furthermore, by Lemma \ref{lem:edge->mi=mj}, the value of $M_i(t)$
is constant on each  of these trees. Finally, all nodes on a
particular tree have either a valid or invalid estimate (the
estimate being valid if and only if $P_r(t)=r$ and
$M_r(t)=u_r(t-1)$ at the root node $r$ of the tree.) For any $z\in
{\cal U}$, we let $G_z(t)$ be the subgraph of $G(t)$ consisting of
those trees at which all nodes have $M_i(t)=z$ and for which the
estimate $z$ on that tree  is invalid. We refer to $G_z(t)$ as the
\emph{invalidity graph} of $z$  at time $t$. In the sequel we will
say that $i$ is in $G_z(t)$, and abuse notation  by writing $i\in
G_z(t)$, to mean that $i$ belongs to the set of nodes  of
$G_z(t)$. The lemmas that follow aim at showing} that the
invalidity graph of the largest estimate eventually becomes empty.
\jred{Loosely speaking,} \jmj{the first \aoc{lemma}} asserts that
\jred{after the $u_i(t)$ have stopped changing,} it takes
essentially two time steps for a maximal estimate $M^*$ to
propagate to a neighboring node.

\begin{lemma} \label{lemma:addition}
\mod{\jnt{Fix some time $t>T'$.} Let $M^*$ be the largest estimate
at \jnt{that time,}  i.e., $M^* = \max_{\jnt{i}} M_i(t)$.
Suppose that $i$ is in $G_{M^*}(t+2)$ but not in $G_{M^*}(t)$. Then
$P_i(t+2)\in G_{M^*}(t)$.}
\end{lemma}
\begin{proof}
The fact that $i\not\in G_{M^*}(t)$ implies that either \aoc{(i)} $M_i(t) \neq
M^*$ or \aoc{(ii)} $M_i(t) = M^*$ and $i$ has a valid estimate at time $t$. In
the latter case, it follows from Lemma
\ref{lem:valid+max-->final} that $i$ also has a valid estimate at
time $t+2$, contradicting the assumption $i\in G_{M^*}(t+2)$. Therefore, \aoc{we can and
will assume that} $M_i(t) <
M^*$. Since $t>T'$, no node ever executes O1. The difference between
$M_i(t)$ and $M_i(t+2)= M^*$ can only result from the execution of
O2 or O4a by $i$ during time slot $t$ or $t+1$.

Node $i$ cannot have executed O4a during slot $t+1$, \aoc{because this would
result in} $P_i(t+2) = i$, and $i$ would have a
valid estimate at time $t+2$, contradicting \aoc{the assumption} $i\in G_M^*(t+2)$.
Similarly, if $i$ executes O4a during slot $t$ it sets $P_i(t+1)
= i$. Unless it executes O2 \aoc{during slot} $t+1$, we have \aoc{again} $P_i(t+2) = i$ contradicting \aoc{
the assumption} $i\in
G_M^*(t+2)$. Therefore, $i$ must have executed O2 during \aoc{either} slot
$t+1$ or slot $t$, and in the latter case it must not
have executed O4a during slot $t+1$.

Let us suppose that $i$ executes O2 during slot $t+1$, and
sets thus $P_i(t+2) =j$ for some $j$ that sent at time $t+1$ a
message with the estimate $M^*=M_i(t+2)$. The rules of the
algorithm imply that {\aoe $M^* = M_j(t+1)$. We can also conclude
that $M_j(t+1)=M_j(t)$, since if this \aoc{were} not true, node $j$ would
have sent out a restart at time $t+1$. Thus $M_j(t) = M^*$. It }
remains to prove that the estimate $M^*$ of $j$ at time $t$ is not
valid. Suppose, to obtain a contradiction, that it is valid. Then it
follows from Lemma \ref{lem:valid+max-->final} that $j$ also has a
valid estimate at time $t+2$, and from the definition of validity
that the estimate of $i$ is also valid at $t+2$, in contradiction
with \aoc{the assumption} $i\in G_{M^*}(t+2)$. {\aoe Thus we have \aoc{established} that }
$P_i(t+2) = j \in G_{M^*}(t)$ if $i$ executes O2 during slot
$t+1$. The same argument applies if $i$ executes O2 during slot
$t$, without executing O4a or O2 during the slot $t+1$, \aoc{using the fact that
in this case} $P_i(t+2) = P_i(t+1)$.
\end{proof}

\jblue{\jred{Loosely speaking,} the \aoc{next} lemma asserts that the
removal of an invalid maximal estimate $M^*$, through the
propagation of restarts, takes place at unit speed.}

\begin{lemma} \label{lemma:roots}
\jblue{Fix some time $t>T'$,} and let $M^*$ be the largest
estimate at \jblue{that time.} Suppose that $i$ is a root
\jnt{(i.e., has zero out-degree)} in the forest $G_{M^*}(t+2)$.
Then, either (i) $i$ is the root of \jnt{an one-element} tree in
\mod{$G_{M^*}(t+2)$} \aoc{consisting} only of $i$, or (ii) $i$ is
\jmj{at least} ``two levels down in $G_{M^*}(t)$, i.e., there
\jnt{exist} nodes $i',i''$ with $(i,i'), (i',i'') \in G_{M^*}(t)$.
\end{lemma}
\begin{proof}
Consider such a node $i$ and assume that (i) does not hold. Then,
\begin{eqnarray} M_i(t) & = & M_i(t+1) = M_i(t+2) = M^* \nonumber \\
P_i(t) & = & P_i(t+1) = P_i(t+2) \label{equality} \end{eqnarray}
\aoc{This is because} otherwise, cases (a) and (b) of Lemma
\ref{lem:speed1/2} imply that $i$ has \aoc{zero in-degree} in
$G_{M^*}(t+2)\subseteq G(t+2)$, in addition to having a zero
out-degree, contradicting our assumption that (i) does not hold.
Moreover, the estimate of $i$ is not valid at $t$, because it
would then also be valid at $t+2$ by Lemma
\ref{lem:valid+max-->final}, in contradiction with $i\in
G_{M^*}(t+2)$. Therefore, $i$ belongs to the forest $G_{M^*}(t)$.
Let $r$ be the root of the connected component to which \aoc{$i$}
belongs. We \aoc{will} prove that $i\neq r$ and $P_i(t) \neq r$,
and thus that (ii) holds.

Since $r\in G_{M^*}(t)$, \aoc{we have} $M_r(t) = M^*$ and
\aoc{thus} $r$ does thus not execute O2 during slot $t$. Moreover,
$r$ is a root and has an invalid estimate, so $P_r(t) \neq r$ and
there is a ``restart" message from $P_r(t)$ to $r$ at time $t$.
Therefore, $r$ executes O4a during slot $t$, setting $P_r(t+1) =r$
and sending ``restart" messages to all its neighbors at time
$t+1$. This implies that $i\neq r$, as we have seen that $P_i(t) =
P_i(t+1)$. Let us now assume, to obtain a contradiction, that
$P_i(t) =r$ and thus \aoc{by Eq. (\ref{equality})}, $P_i(t+2) =
P_i(t+1) =r$. In that case, we have just seen that there is at
time $t+1$ a ``restart" message from $r=P_i(t+1)\neq i$ to $i$, so
$i$ executes O4a during slot $t+1$ and sets $P_i(t+2) =i$. This
however contradicts the fact $P_i(t+2) = P_i(t+1)$. Therefore,
{\aoe $r \neq i$ and $r \neq P_i(t)$,} i.e.,  $i$ is ``at least
two levels down" in $G_{M^*}(t)$.
\end{proof}

\ao{Let the {\em depth} of a tree be the largest distance between a leaf
of the tree and the root; the depth of a forest is the largest depth
of any tree in the forest. We will use $g(\cdot)$ to denote depth. }
\jblue{The following lemma uses the previous two lemmas to
assert that a {\aoe forest }carrying an invalid maximal estimate
has its depth decrease by at least one over a time interval
of length two.}

\begin{lemma} \jblue{Fix some
time $t > T'$,} and let $M^*$ be the largest estimate value at
\jblue{that time.} If $g(G_{M^*}(t+2))
> 0$, then \aoc{$g( G_{M^*}(t+2)) \leq g(G_{M^*}(t))-1.$}
\label{lemma:depth_decrease}
\end{lemma}

\begin{proof}
\jblue{Suppose that $g(G_{M^*}(t+2)) > 0$. Let us fix a leaf $i$
and a root $j$ in the forest $G_{M^*}(t+2)$ such that the length
of the path from $i$ to $j$ is equal to the depth of
$G_{M^*}(t+2)$.} \mod{Let $i'$ be the \jblue{{\aoe \aoc{single
neighbor}} of node} $i$ in $G_{M^*}(t+2)$. \jblue{We \jmj{first}
claim} that every edge $(k,k')$ on the path from $i'$ to $j$ in
$G_{M^*}(t+2)$ was also present in $G_{M^*}(t)$.} { \aoe Indeed,
by Lemma \ref{lem:speed1/2}, the appearance of a new edge $(k,k')$
at time $t+1$ or $t+2$ implies that node $k$ has in-degree $0$ in
$G(t+2)$, which contradicts $k$ being an intermediate node on the
path from $i$ to $j$ in $G_{M^*}(t+2)$. The same argument
establishes that $M_k(t)=M_k(t+1)=M^*$. Finally, the estimate of
$k$ at time $t$ is \aoc{invalid}, for if it were \aoc{valid}, it
would still be \aoc{valid} at time $t+2$ by Lemma
\ref{lem:valid+max-->final}, so $i$ \aoc{would also} have a
\aoc{valid} estimate at time $t+2$, which is false by assumption.
Thus we have just established \aoc{that both} the node $k$ and its
edge $(k,k')$ at time $t+2$ belong to $G_{M^*}(t)$. }

\aoc{Thus the graph $G_{M^*}(t)$ includes a path from $i'$ to $j$ of length
$g(G_{M^*}(t+2))-1$.} Moreover, by Lemma \ref{lemma:roots}, we know
that at time $t$ \jblue{some} edges $(j,j')$ and $(j',j'')$ were
present in $G_{M^*}(t)$, so the \jred{path length} from $i'$ to
\jmj{$j''$ is at least $g( G_{M^*}(t+2)) + 1$}. This proves \aoc{
that $g( G_{M^*}(t)) \geq g(G_{M^*}(t+2))+1$ and the
lemma.}
\end{proof}

The following lemma {\aoe analyzes the remaining case of invalidity
graphs with zero depth.} \aoc{It shows that the invalidity graph will be empty
two steps after its depth reaches zero.}

\begin{lemma}\label{lem:invalid_length0}
Fix some time $t > T'$, and let $M^*$ be the largest estimate value
at \jblue{that time.}  If $G_{M^*}(t+2)$ is not empty, then $g(
G_{M^*}(t+1))
>0 $ or $g( G_{M^*}(t))  >0 $.
\end{lemma}

\begin{proof}
Let us take a node $i\in G_{M^*}(t+2)$ and let $j = P_i(t+2)$. It follows from the
definition of a valid estimate that $j \neq i$. This
implies that $i$ did not execute O4a (or O1) during slot $t+1$.
We treat two cases separately:

(i) \emph{Node $i$ did not execute O2 during slot $t+1$.}
In this case, $P_i(t+1) = P_i(t+2) = j$ and $M_j(t+1) = M_j(t+2) =
M^*$. Besides, there is no ``restart" message from $j=P_i(t+1)$ to
$i$ at time $t+1$, for otherwise $i$ would have executed O4a
during slot $t+1$, which we know it did not. Therefore,  $(i,j)\in
E(t {\aoe + 1})$ by definition of $G(t)$, and $M_j(t+1) = M_i(t+1)
= M^*$ by Lemma \ref{lem:edge->mi=mj}. Moreover, neither $i$ nor
$j$ have a valid estimate, for otherwise Lemma
\ref{lem:valid+max-->final} would imply that they both hold the
same valid estimate at $t+2$, in contradiction with $i\in
G_{M^*}(t+2)$. So the edge $(i,j)$ is present in $G_{M^*}(t+1)$,
which has thus a positive minimal depth.\\
(ii) \emph{Node $i$ did execute O2 during slot $t+1$:}
In that case, there was a message \aoc{with} the value $M_i(t+2) = M^*$
from $j$ to $i$ at time $t+1$, which implies that $M_j(t+1) =
M_j(t+2) = M^*$. This implies that $j$ did not execute operation O2
during slot $t$. Moreover, node $j$ did not have a valid
estimate at time $t+1$. Otherwise, {\aoe part (a) }of Lemma
\ref{lem:valid+max-->final} implies that \aoc{$j$ has} a valid
estimate at time $t+2$, {\aoe and part (b) of the same lemma
implies} there \aoc{was not} a ``restart" message {\aoe from} $j$ at
$t+2$, so that $(i,j)\in E(t+2)$. This would in turn imply that $i$
has a valid estimate at time $t+2$, contradicting $i\in
G_{M^*}(t+2)$. {\aoe To summarize,} $j$ has an invalid estimate
$M^*$ at time $t+1$ and did not execute O2 during slot $t$.
{\aoe We now simply observe that the argument of case (i) applies}
to $j$ at time $t+1$.
\end{proof}

\jmj{The next lemma asserts that the largest invalid estimates are
eventually purged, and thus that eventually, all remaining largest
estimates are valid.}

\begin{lemma}\label{prop:no_unreliable}
Fix some time $t > T'$, and let $M^*$ be the largest estimate value
at \jblue{that time.} Eventually $G_{M^*}(t)$ is empty.
\end{lemma}
\begin{proof}
{\aoe Lemma \ref{lemma:depth_decrease} implies there is a time
$t'>T'$ after which $g(G_{M^*}(t'')) =0$ for all $t''>t'$. Lemma
\ref{lem:invalid_length0} then implies} that $G_{M^*}(t)$ is empty
for all $t>t'+2$.
\end{proof}

\jnt{We are now ready for the {\aoe proof of the main theorem}.}

\begin{proof}[Proof of Theorem \ref{thm:maxtracking}.]
\jblue{Let $\oM=\max_i u_i(T')$.} \jmj{It follows from the
definition of a valid estimate that any node holding an estimate
$M_i(t)> \oM$ \aoc{at time $t \geq T'$} has an invalid estimate.} Applying Lemma
\ref{prop:no_unreliable} repeatedly shows the existence of some time
$\oT\jblue{\geq T'}$ such that when $t \geq \oT$, \jmj{no node has
an estimate larger than $\oM$, and every node having an estimate
$\oM$ has a valid estimate.}

{\aoe We will assume that the time $t$ in
every statement we make below satisfies $t \geq \overline{T}$.  Define} $Z(t)$ as the set of nodes
having the estimate $\oM$ at time $t$. Every node in $Z(t)$ holds a
valid estimate, and $Z(t)$ is never empty because Lemma
\ref{lem:m_geq_y} implies that $M_i(t)\geq \oM$ for every $i$ with
$u_i = \oM$. Moreover, it follows from Lemma
\ref{lem:valid+max-->final} and the definition of validity  that any
node belonging to some $Z(t)$ will forever afterwards \aoc{maintain $M_i(t) = \oM$ and will} satisfy
 the conclusion of  Theorem
\ref{thm:maxtracking}.

{\aoe We conclude the proof by arguing that eventually every node is
in $Z(t)$. In particular, we will argue that}  a node $i$ adjacent
to a node $j\in Z(t)$ necessarily belongs to $Z(t+2)$. {\aoe Indeed,
} it follows from Lemma \ref{lem:valid+max-->final} that node $j$
sends to $i$ a message \aoc{with} the estimate $\oM$ at time $t+1$. If
$i\in Z(t+1)$, then $i \in Z(t+2)$; else $M_i(t+1)< \oM$, $i$
executes O2 during slot $t+1$, and sets $M_i(t+2) = \oM$, so indeed
$i\in Z(t+2)$. \end{proof}

%% file: complexity.tex
\def\aoc#1{#1}

\section{The Time to Termination}\label{sec:complexity}

The aim of this appendix is to prove an $O(n^2 K^2)$ upper bound
on the time to termination of the interval-averaging algorithm
from Section \ref{sec:comput_average}. The validity of Theorems
\ref{averagingtime} and \ref{approximatetime} follows readily
because finite-memory algorithms described in Section
\ref{sec:reduction_to_avg} all consist \aoc{of} having the nodes
run a fixed number of interval-averaging algorithms.

We will use the notation $M'(t)$ to denote the largest estimate
\jh{held by any node} at time $t$ or in the $n$ time steps
preceding it: \[ M'(t) = \max_{     \begin{footnotesize}
\begin{array}{c}
                                    i=1,\ldots,n \\
                                    k=t,t-1,\ldots,t-n \\
                                  \end{array} \end{footnotesize}
} ~~~M_i(k).\] For $M'(t)$ to be well defined, we will adopt the
convention that for all negative times $k$, $M_i(k)=u_i(0)$.

\aoc{\begin{lemma}\label{lem:M'monotonous} In the course of the
execution of the interval-averaging algorithm, $M'(t)$ never
increases.
\end{lemma} }
\begin{proof}
Fix a time $t$. We will argue that \begin{equation} \label{eq:it}
M_i(t+1) \leq M'(t) \end{equation} for each $i$. This clearly
implies $M'(t+1) \leq M'(t)$.

If $M_i(t+1) \leq M_i(t)$, then Eq. (\ref{eq:it}) is obvious. We
can thus suppose that $M_i(t+1) > M_i(t)$. There are only three
boxes in Figure \ref{fig:maxtracking} which result in a change
between $M_i(t)$ and $M_i(t+1)$. These are $O2$, $O1$, and $O4a$.
We can rule out the possibility that node $i$ executes $O2$, since
that merely sets $M_i(t+1)$ to some $M_j(t)$, and thus cannot
result in $M_i(t+1) > M'(t)$.

Consider, then, the possibility that node $i$ executes $O1$ or
$O4a$, and as a consequence $M_i(t+1)=u_i(t)$. If $u_i(t) \leq
u_i(t-1)$, then we are finished because \[ M_i(t+1) = u_i(t) \leq
u_i(t-1) \leq M_i(t), \] which contradicts the assumption
$M_i(t+1)>M_i(t)$. Note that the last step of the above chain of
inequalities used Lemma \ref{lem:m_geq_y}.

Thus we can assume that $u_i(t)>u_i(t-1)$. In this case, $i$ must
have fulfilled acceptance from some node $j$ during slot $t-1$.
\jh{Let $\hat{t}$ the time when node $j$ received the
corresponding request message from node $i$. The rules of our
algorithm imply that $u_i(\hat t) = u_i(t-1)$, and that the
quantity $w$ sent by $j$ to $i$ in response to the request is no
greater than $\frac{1}{2}(u_j(\hat t)-u_i(\hat t))$. This implies
that $u_i(t) = u_i(t-1) +w < u_j(\hat t)$. \comjh{I've add some
more explanations on why $u_i(t)<u_j(\hat t)$}}

Crucially,  we have that $\hat{t} \in [t-1-n,t-1]$, since at most
$n+1$ time steps pass between the time node $j$ receives the
request message it will accept and the time when node $i$ fulfills
$j$'s acceptance. So \[ M_i(t+1) = u_i(t) < u_j(\hat{t}) \leq
M_j(\hat{t}+1) \leq M'(t).\] We have thus showed that $M_i(t+1)
\leq M'(t)$ in every possible case, \jh{which implies that
$M'(t+1)\leq M'(t)$.}
\end{proof}

\begin{lemma}\label{lem:conv_time_minmax}
\ora{Consider the maximum tracking algorithm.} If each
\aoc{$u_i(t)$} is constant for $\aoc{t \in } [t_0,t_0+4n]$, then
\aoc{at least one of the following two statements is true}:
\texitem{(a)} $M'(t_0+3n) \ora{<} M'(t_0)$.

\texitem{(b)} $M_i(t_0+4n) = \max_j u_j(t_0)$ for every $i$.
\end{lemma}
\begin{proof}
\jh{Suppose first that no node holds an estimate equal to
$M'(t_0)$ at some time between $t_0+2n$ and $t_0+3n$. Then it
follows from the definition of $M'(t)$ and its monotonicity
(\ref{lem:M'monotonous}) that condition (a) holds. Suppose now
that some node holds an estimate equal to $M'(t_0)$ at some time
between $t_0+2n$ and $t_0+3n$. The definition of $M'(t)$ and the
monotonicity of $\max_i M_i(t)$ when all $u_i$ are constant (Lemma
\ref{lem:nondecreasing}) imply that $M'(t_0) = \max_i
M_i(t)$ for all $t\in [t_0,t_0+3n]$. It follows from repeated
application of Lemmas \ref{lemma:depth_decrease} and
\ref{lem:invalid_length0} (similarly to what is done in the proof
of Proposition \ref{prop:no_unreliable}) that every estimate
$M'(t_0)$ at time $t_0 + 2n$ is valid, which by definition implies
the existence of at least one node $i$ with $u_i(t_0) = M'(t_0)$.
Besides, since $M_i(t)\geq u_i(t)$ holds for all $i$ and $t$ by
Lemma \ref{lem:m_geq_y} and since we know that $M'(t_0) = \max_i
M_i(t)$ for all $t\in [t_0,t_0+3n]$, we have $M'(t_0) = \max_i
u_i(t_0)$. As described \ora{in the} proof of Theorem
\ref{thm:maxtracking}, this implies that after at most $2n$ more
time steps, $M_i(t)=M'(t_0)$ holds for every $i$, \aoc{and so (b)
holds.}}
\end{proof}

The next lemma \aoc{upper bounds} the largest time before
\aoc{some} request is accepted or \aoc{some} outdated estimate is
purged from the system. Recall that $\bar x = (\sum_i^nx_i)/n$.

\begin{lemma}\label{lem:bound_interoperation_time}
\ora{Consider the interval-averaging} algorithm described in
Sections \ref{sec:max_tracking} and \ref{sec:comput_average}. For
any $t_0$, at least one of the following is true: \texitem{(a)}
Some node accepts a request at \ora{some slot $t\in
[t_0,t_0+8n-1]$.} \texitem{(b)} We have $M'(t+1) < M'(t)$ for some
$t\in [t_0+1,t_0+3n]$. \texitem{(c)} All $u_i$ remain forever
\aoc{constant} after time $t_0+n$, with $u_i \in \{\floor{\bar x}
, \ceil{\bar x}\}$, and all $M_i(t)$ remain \aoc{forever} constant
after time $t_0 + 5n$, with $M_i = \ceil{\bar x}$.
\end{lemma}
\begin{proof}
Suppose that condition (a) does not hold, i.e., that no node
accepts a request between $t_0$ and $t_0 + 8 n$. \ora{Since an
acceptance message needs to travel through} at most $n-2$
intermediate nodes before reaching the \ora{originator} of the
request (cf.\ Lemma \ref{lem:request_terminate}), we conclude that the system does
not contain any \ora{acceptance messages}  after \ora{time}
$t_0+n$. As a result, no node modifies its \ora{value} $u_i(t)$
between \ora{times} $t_0+n$ and $t_0 + 8 n$. Abusing notation \aoc{slightly}, we
\aoc{will} call these values $u_i$.

It follows from Lemma \ref{lem:conv_time_minmax} that either
condition (b) holds, or that there is a time $\tilde{t} \leq
t_0+5n$ at which $M_i(\tilde{t}) = \max_j u_j$ for every $i$. Some
requests may have been emitted in the interval $[t_0,t_0+5n]$.
\aoc{Since we are assuming that} condition (a) does not hold,
\ora{these requests must have all been} denied. It follows from
Lemma \ref{lem:request_terminate} that none of {these requests} is
present by time $t_0+ 7n$. Moreover, by the rules of our
algorithm, once $M_i$ \ora{becomes equal to} $\max_j u_j$ for
every node $i$, every node with $u_i \leq \max_j u_j -2$ \ora{will
keep} emitting requests. \ora{Using an argument similar to the one
at the end of the proof of}  Proposition
\ref{prop:final_proof_average}, if such requests are sent, at
least one must be accepted within $n$ time steps, that is, no
later than \ora{time} $t_0 +8n$. \aoc{Since by assumption this has
not happened, no such requests could have been sent, implying
that} $u_i \geq \max_j u_j -1$ for every $i$. \aoc{Moreover, this
implies that no request messages/acceptance are ever sent after
time $t_0 + 7n$, so that $u_i$ never change. It is easy to see
that the $M_i$ never change as well, so that condition condition
(c) is satisfied.}
\end{proof}

We can now give an upper bound on the \ora{time until our
algorithm terminates.}

\begin{theorem}
The interval-averaging algorithm described in \aoc{Section
\ref{sec:comput_average}} \ora{terminates after} at most
\aoc{O($n^2K^2)$} time steps.
\end{theorem}
\begin{proof}
Consider the function $V(t) = \sum_{i=1}^n \left(\hat u_i(t) -
\bar u\right)^2$, where $\bar u$ \ora{is the average of the $x_i$,
which is also the average of the $\hat u_i$, and where the $\hat
u_i(t)$ are} as defined before Lemma
\ref{prop:finite_request_acceptances}. Since $\hat u_i(t)\in
\{0,1,\dots,K\}$ for all $i$, one can verify that $V(0) \leq
\frac{1}{4}n K^2$.  Moreover, as explained in the proof of Lemma
\ref{prop:finite_request_acceptances}, $V(t)$ is non-increasing,
and decreases by at least 2 with every request acceptance.
Therefore, \ora{a total of} at most $\frac{1}{8}nK^2$ requests can
be accepted. Furthermore, \aoc{we showed that $M'(t)$ is
non-increasing, and since $M'(t)$ always belongs to
$\{0,1,\dots,K\}$, it can strictly decrease at most $K$ times.} It
follows then from Lemma \ref{lem:bound_interoperation_time} that
condition (c) must hold after at most $\frac{1}{8}nK^2\cdot 8n +
K\cdot 3n \jh{+5}$ time steps.

\ora{Recall that} in parallel with the maximum-tracking and
averaging algorithm, we \ora{also} run a minimum tracking
algorithm. In the previous paragraph, we demonstrated that
condition (c) of Lemma \ref{lem:bound_interoperation_time} holds,
i.e. $u_i$ remain fixed forever, after $n^2 K^2 + K\cdot 3n$ time
steps. A similar argument to Lemma \ref{lem:conv_time_minmax}
implies that the minimum algorithm will reach a fixed point after
an additional \jh{$(3K+4)n$} steps.

Putting it all together, the algorithm reaches a fixed point after
\jh{$n^2 K^2 + (6K+4)\cdot n$} steps.
\end{proof}

\bigskip

\begin{figure}
\centering
\includegraphics[scale= .45]{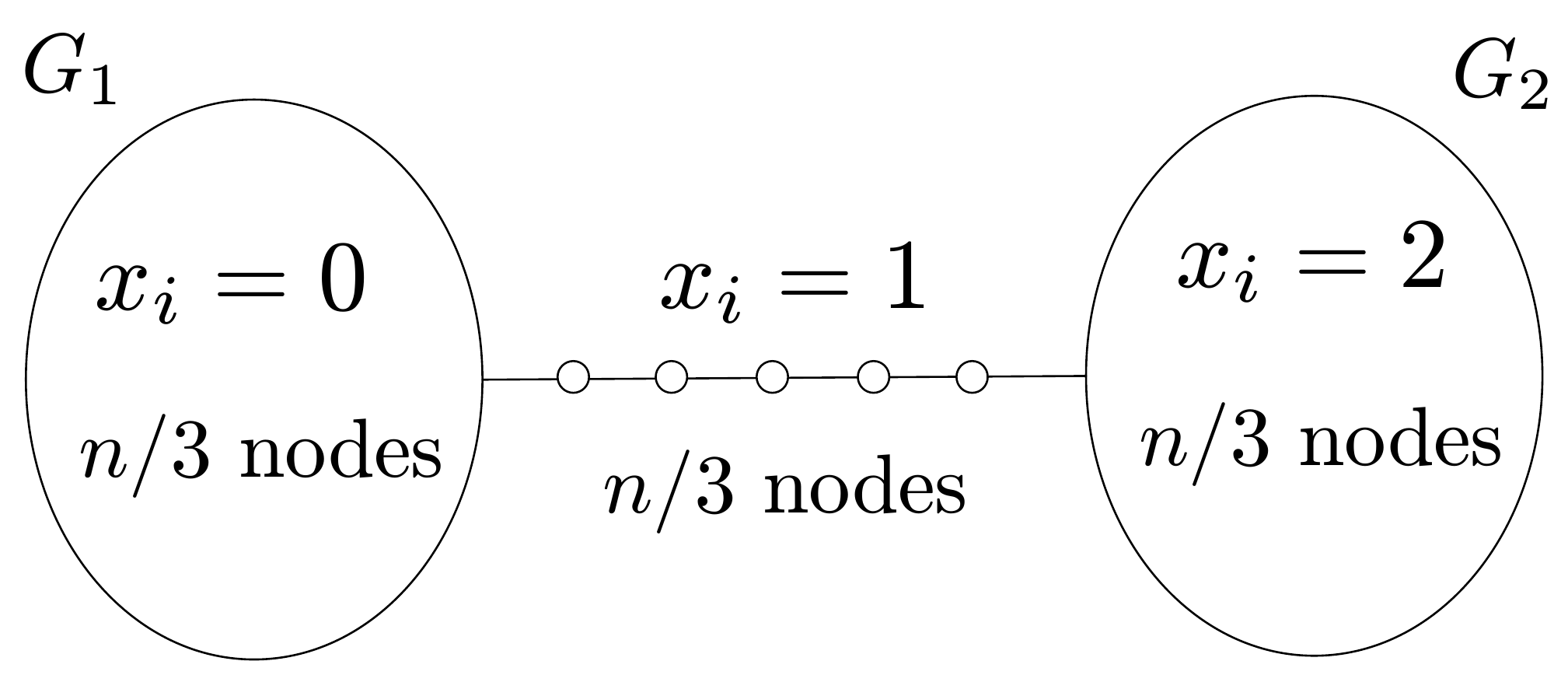}
\caption{A class of networks and initial
conditions for which our algorithm takes $\aoc{\Theta}(n^2)$ time steps to
reach its final state.
}\label{fig:lower_complexity_bound}
\end{figure}

\aoc{We note that there are cases where the running time of
interval averaging is quadratic in $n$. For example,} consider the
network in Figure \ref{fig:lower_complexity_bound}, consisting of
two arbitrary \aoc{connected} graphs $G_1,G_2$ with $n/3$ nodes
each, connected by a \aoc{line} graph of $n/3$ nodes. Suppose that
$K=2$, and that $x_i=0$ if $i\in G_1$, $x_i=2$ if $i\in G_2$, and
$x_i =1$ otherwise. \aoc{The algorithm will have} the nodes of
$G_1$ with $u_i = 0$ send requests to nodes $j$ in $G_2$ with
$u_j=2$, and each successful request will result in the pair of
nodes changing their values, $u_i$ and $u_j$, to 1. The system
will reach its final state after $n/3$ such successful requests.
Observe now that each successful request must cross the line
graph, which takes $n/3$ time steps in each direction. Moreover,
since nodes cannot simultaneously treat multiple requests, once a
request begins crossing the line graph, all other requests are
denied until the response to the first request reaches $G_1$,
which takes at least $2n/3$ time steps. Therefore, in this
example, it takes at least $2n^2/\aoc{9}$ time steps until the
algorithm terminates.

